\documentclass{amsart}
\usepackage{amssymb,latexsym,amsmath,amsfonts,hhline}
\usepackage{pdfsync}
\usepackage{color}

\newcommand{\cal}{\mathcal}

\makeatletter
\renewcommand{\subsection}{\@startsection{subsection}{2}{0mm}{-2mm}{-2mm}{\bf\normalsize}}

\def\sbsnt#1{\subsection{#1}}
\makeatother

\newtheorem{formula}{}[section]
\newtheorem{definition}[formula]{Definition}
\newtheorem{corollary}[formula]{Corollary}
\newtheorem{remark}[formula]{Remark}
\newtheorem{lemma}[formula]{Lemma}
\newtheorem{theorem}[formula]{Theorem}

\def\thrm{\begin{theorem}}
\def\thrml#1{\begin{theorem}\label{#1}}
\def\ethrm{\end{theorem}}
\def\rmrk{\begin{remark}}
\def\rmrkl#1{\begin{remark}\label{#1}}
\def\ermrk{\end{remark}}
\def\dfntn{\begin{definition}}
\def\dfntnl#1{\begin{definition}\label{#1}}
\def\edfntn{\end{definition}}
\def\nmrt{\begin{enumerate}}
\def\enmrt{\end{enumerate}}
\def\tm#1{\item[{\rm (#1)}]}
\def\qtn{\begin{equation}}
\def\qtnl#1{\begin{equation}\label{#1}}
\def\eqtn{\end{equation}}
\def\lmm{\begin{lemma}}
\def\lmml#1{\begin{lemma}\label{#1}}
\def\elmm{\end{lemma}}
\def\crllr{\begin{corollary}}
\def\crllrl#1{\begin{corollary}\label{#1}}
\def\ecrllr{\end{corollary}}
\def\css{\begin{cases}}
\def\ecss{\end{cases}}

\def\proof{\noindent{\bf Proof}.\ }

\def\cA{{\cal A}}
\def\cB{{\cal B}}

\def\cS{{\cal S}}

\def\cX{{\cal X}}

\def\mF{{\mathbb F}}

\def\mZ{{\mathbb Z}}

\DeclareMathOperator{\aut}{Aut}

\DeclareMathOperator{\cyc}{Cyc}

\DeclareMathOperator{\iso}{Iso}

\DeclareMathOperator{\orb}{Orb}

\DeclareMathOperator{\rad}{rad}
\DeclareMathOperator{\rk}{rk}

\DeclareMathOperator{\Span}{Span}

\DeclareMathOperator{\sym}{Sym}
\DeclareMathOperator{\tr}{tr}

\def\bull{\hfill\vrule height .9ex width .8ex depth -.1ex\medskip}

\def\mmod#1#2#3{#1=#2\ (\text{\rm mod}\hspace{2pt}#3)}
\def\nmmod#1#2#3{#1\ne #2\ (\text{\rm mod}\hspace{2pt}#3)}
\def\qaq{\quad\text{and}\quad}
\def\qoq{\quad\text{or}\quad}

\newcommand{\grp}[1]{\langle {#1}\rangle}
\newcommand{\und}[1]{{\underline{#1}}}
\def\twoe{\approx_{\scriptscriptstyle 2}}                           

\begin{document}
\title{On Schur 2-groups}
\author{Mikhail Muzychuk}
\address{Netanya Academic College, Netanya, Israel}
\email{muzy@netanya.ac.il}
\author{Ilya Ponomarenko}
\address{Steklov Institute of Mathematics at St. Petersburg, Russia}
\email{inp@pdmi.ras.ru}
\thanks{The work of the second author was partially supported by the RFBR Grant 14-01-00156 А}
\date{}

\begin{abstract}
A finite group $G$ is called a Schur group, if any Schur ring over~$G$ is the transitivity module of a point stabilizer
in a subgroup of $\sym(G)$ that contains all right translations. We complete a classification of abelian $2$-groups by proving that the 
group $\mZ_2\times\mZ_{2^n}$ is Schur. We also prove that any non-abelian Schur $2$-group of order larger than $32$ is dihedral 
(the Schur $2$-groups of smaller orders are known). Finally, in the dihedral case, we study Schur rings of rank at most $5$, and show that
the unique obstacle here is a hypothetical  S-ring of rank $5$ associated with a divisible difference set.
\end{abstract}

\maketitle

\section{Introduction}
Following R.P\"oschel \cite{Poe74}, a finite group $G$ is called a {\it Schur group}, if any S-ring over~$G$ is the transitivity module 
of a point stabilizer in a subgroup of $\sym(G)$ that contains all right translations (for the exact definitions, we refer to 
Section~\ref{120813a}). He  proved there that if $p\ge 5$ is a prime, then a finite $p$-group is Schur if and only if 
it is cyclic. For $p=2$ or $3$, a cyclic $p$-group is still Schur, but  P\"oschel's
theorem is  not true: a straightforward computation shows that an elementary abelian group
of order $4$ or $9$ is Schur. In this paper, we are interested in Schur 2-groups.\medskip

Recently in \cite{EKP14}, it was proved that every finite abelian Schur group belongs to one 
of several explicitly given families. In particular, from Lemma~5.1 of that paper it 
follows that all abelian Schur $2$-groups are known except for the groups 
$\mZ_{2^{}}\times \mZ_{2^n}$, where $n\ge 5$. We prove that all these groups
are Schur (Theorem~\ref{090514b}). As a by-product we can complete the classification of 
abelian Schur $2$-groups.

\thrml{170814a}
An abelian $2$-group $G$ is Schur if and only if $G$ is cyclic, or elementary abelian of order at most $32$, or
is isomorphic to $\mZ_{2^{}}\times \mZ_{2^n}$ for some~$n\ge 1$.
\ethrm

Non-abelian Schur groups have been studied in~\cite{PV} where it was proved that they are
metabelian. In particular, from Theorem~4.2 of that paper 
it follows that non-abelian Schur $2$-groups are known except for dihedral groups and groups
\qtnl{170914a}
M_{2^n}=\grp{a,b:\ a^{2^{n-1}}=b^2=1,\ bab=a^{1+2^{n-2}}},
\eqtn
where $n\ge 4$. In this paper we prove that the latter groups are not Schur (Theorem~\ref{071113a}).
As a by-product we obtain the following statement.

\thrml{170814b}
A non-abelian Schur $2$-group of order at least $32$ is dihedral.
\ethrm

We do not know whether or not a dihedral $2$-group of order more than $32$ is Schur. A  standard technique based 
on Wielandt's paper~\cite{W49} enables us to describe S-rings of rank at most~$5$ as follows (see  Subsection~\ref{231014w} for a connection between S-rings and divisible difference sets). 

\thrml{170814c}
Let be $\cA$ an S-ring  over  a dihedral $2$-group. Suppose that
$\rk(\cA)\le 5$. Then one of the following statements is true:
\nmrt
\tm{1} $\cA$ is isomorphic to an S-ring over $\mZ_{2^{}}\times \mZ_{2^n}$,
\tm{2} $\cA$ is a proper dot or wreath product,
\tm{3} $\rk(\cA)=5$ and $\cA$ is associated with a divisible difference set in $\mZ_{2^n}$.
\enmrt
\ethrm

S-rings in statement (1) of this theorem are schurian by Theorem~\ref{170814a}. By 
induction, this implies that all S-rings in statement~(2) are also schurian. Thus, by
Theorem~\ref{081014w} we obtain the following corollary.

\crllrl{031014a}
Under the hypothesis of Theorem~\ref{170814c}, the S-ring $\cA$ is 
not schurian only if $\cA$ is associated with a divisible difference set in a cyclic $2$-group.
\ecrllr

In fact, we do not know whether there exists a non-trivial divisible difference set in a cyclic $2$-group that produces an S-ring $\cA$ 
in part~(3) of Theorem~\ref{170814c}. If such a set does exist, then the corresponding dihedral $2$-group is not Schur. On the other hand, 
in Subsection~\ref{231014w}, we show that using a relative difference set (which is a special case of a divisible one), one can construct 
an S-ring of rank~$6$ (over a dihedral $2$-group). These difference sets, and, therefore, S-rings, do exist, but are relatively rare. The only known example is the classical $(q+1,2,q,(q-1)/2)$-difference set
where $q$ is a Mersenne prime. Thus, the question whether a dihedral $2$-group is Schur, 
remains open.\medskip

The paper consists of fourteen sections. In Sections~\ref{120813a}, \ref{111014a}
and \ref{231014v}, we give a background of S-rings, Cayley schemes\footnote{Here, we assume
some knowledge of association scheme theory.} and cite
some basic facts on S-rings over cyclic and dihedral groups. In Sections~\ref{231014a1}--
\ref{231014a3}, we develop a theory of S-rings over $\mZ_2\times\mZ_{2^n}$ that is
culminated in Theorem~\ref{090514b} stating that this group is Schur.
In Section~\ref{231014a4}, we show that the group $M_{2^n}$ is not Schur 
for all $n\ge 4$ (Theorem~\ref{071113a}).
In Sections~\ref{231014a6}--\ref{231014a7}, we study S-rings over a dihedral $2$-group:
here, we start with constructions based on cyclic divisible difference sets, and then
complete the proof of Theorem~\ref{170814c}.\medskip

{\bf Notation.}
As usual, by $\mZ$ we denote the ring of rational integers.

The identity of a group $D$ is denoted by $e$; the set of non-identity elements in $D$ is denoted by  $D^\#$.

Let $X\subseteq D$. The subgroup of $D$ generated by $X$ is denoted by $\grp{X}$; 
we also set $\rad(X)=\{g\in D:\ gX=Xg=X\}$. The element $\sum_{x\in X}x$ of the group 
ring $\mZ D$ is denoted by $\und{X}$. The set $X$ is called {\it regular} if the order
$|x|$ of an element $x\in X$ does not depend on the choice of~$x$.

For a group $H\trianglelefteq D$, the quotient epimorphism from $D$ onto $D/H$ is denoted by $\pi_{D/H}$.

The group of all permutations of $D$ is denoted by $\sym(D)$. The set of orbits of a group 
$G\le\sym(D)$ is denoted by $\orb(G)=\orb(G,D)$. We write $G\twoe G'$ if the groups $G,G'\le\sym(D)$
are {\it $2$-equivalent}, i.e. have the same orbits in the coordinate-wise action on~$D\times D$.

Given two subgroups $L\trianglelefteq U\leq D$,  the quotient group $U/L$ is called the {\it section} of $D$. For a 
set $\Delta\subseteq\sym(D)$ and a section $S = U/L$ of  $D$ we set
$$
\Delta^S=\{f^S:\ f\in \Delta,\ S^f=S\},
$$
where $S^f = S$ means that $f$ permutes the right $L$-cosets in $U$ and $f^S$ denotes the bijection of 
$S$ induced by~$f$.

The cyclic group of order $n$ is denoted by  $\mZ_n$.

\section{A background on S-ring theory}\label{120813a}

In what follows,  we use the notation and terminology of~\cite{EP12}.\medskip

Let $D$ be a finite group. A subring~$\cA$ of the group ring~$\mZ D$ is called a {\it Schur 
ring} ({\it S-ring}, for short) over~$D$ if there exists a partition $\cS=\cS(\cA)$ of~$D$ 
such that
\nmrt
\tm{S1} $\{e\}\in\cS$,
\tm{S2} $X\in\cS\ \Rightarrow\ X^{-1}\in\cS$,
\tm{S3} $\cA=\Span\{\und{X}:\ X\in\cS\}$.
\enmrt
When $\cS=\orb(K,D)$ where $K\le\aut(D)$, the S-ring $\cA$ is called {\it cyclotomic} and
denoted by $\cyc(K,D)$. A group isomorphism $f:D\to D'$ is called a {\it Cayley isomorphism} from an S-ring $\cA$ over 
$D$ to an S-ring $\cA'$ over $D'$ if $\cS(\cA)^f=\cS(\cA')$.\medskip

It follows from (S3) that given $X,Y\in\cS(\cA)$ there exist non-negative integers
$c_{X Y}^Z$, $Z\in\cS(\cA)$, such that
$$
\und{X}\,\und{Y}=\sum_{Z\in\cS(\cA)}c_{X Y}^Z\und{Z}.
$$
One can see that  $c_{XY}^Z$ equals the number of different representations $z=xy$ with 
$(x,y)\in X\times Y$ for a fixed (and hence for all) $z\in Z$. It is a well-known fact that 
$$
c_{Y^{-1}X^{-1}}^{Z^{-1}}=c_{X^{}Y^{}}^{Z^{}}\qaq |Z|c_{XY}^{Z^{-1}}=|X|c_{YZ}^{X^{-1}}=|Y|c_{ZX}^{Y^{-1}}
$$ 
for all $X,Y,Z$. A ring isomorphism $\varphi:\cA\to \cA'$ is said to be {\it algebraic} if for any $X\in\cS(\cA)$
there exists $X'\in\cS(\cA')$ such that $\varphi(\und{X})=\und{X'}$.\medskip

The classes of the partition $\cS$ and the number $\rk(\cA)=|\cS|$ are called the {\it basic 
sets} and the {\it rank} of the S-ring~$\cA$, respectively. Any union of basic sets
is called an {\it $\cA$-subset of~$D$} or {\it $\cA$-set}. The set of all of them
is closed with respect to taking inverse and product. Given an $\cA$-set $X$, we denote by $\cA_X$
the submodule of~$\cA$ spanned by  the elements $\und{Y}$, where $Y$ belongs to the set
$$
\cS(\cA)_X=\{Y\in\cS(\cA):\ Y\subseteq X\}.
$$

Any subgroup of $D$ that is an $\cA$-set, is called an {\it $\cA$-subgroup} of~$D$ or 
{\it $\cA$-group}. With each $\cA$-set $X$, one can  naturally associate two $\cA$-groups, 
namely $\grp{X}$ and  $\rad(X)$ (see Notation). The following useful lemma was proved 
in \cite[p.21]{EvdP09}.\footnote{Apparently, for the first time this lemma was proved in~\cite[Proposition~4.5]{M87}.}

\lmml{090608a}
Let $\cA$ be an S-ring over a group $D$ and $H\le D$ an $\cA$-group. Then
given $X\in\cS(\cA)$, the cardinality of the set $X\cap xH$ does not depend on $x\in X$.
\elmm

A section $S=U/L$ of the group $D$ is  called an  {\it $\cA$-section}, if both  $U$ and 
$L$ are $\cA$-groups. In this case,  the module
$$
\cA_S=\Span \{\pi_S(X):\ X\in\cS(\cA)_U\}
$$
is an S-ring over the group~$S$, the basic sets of which are exactly the sets $\pi_S(X)$ from the 
right-hand side of the formula.\medskip

The S-ring $\cA$ is called {\it primitive} if the only $\cA$-groups are $e$ and $D$,
otherwise this ring is called  {\it imprimitive}. One can see that if $H$
is a minimal $\cA$-group, then the S-ring $\cA_H$ is primitive. The classical results
on primitive S-rings over abelian and dihedral groups were obtained in papers~\cite{S33,K37,W49}. 
A careful analysis of the proofs shows that the schurity assumption there
was superfluous. Therefore, in the first part of the following statement, we formulate the
corresponding results in slightly more general form.

\thrml{030414a}
Let $D$ be a $2$-group which is cyclic, dihedral, or isomorphic to the group $\mZ_2\times\mZ_{2^n}$. Then any primitive S-ring over $D$
is of rank~$2$. In particular, if $\cA$ is an S-ring over $D$ and $H\le D$ is a minimal $\cA$-group,
then $H^\#$ is a basic set of~$\cA$. 
\ethrm

Let $S=U/L$ be an $\cA$-section of the group $D$. The S-ring $\cA$ is called the {\it generalized $S$-wreath product} 
\footnote{In \cite{LM98}, the term {\it wedge product} was used.} if the group $L$ is
normal in~$D$ and $L\le\rad(X)$ for all basic sets $X$ outside~$U$; in this case we write
\qtnl{050813a}
\cA=\cA_U\wr_S\cA_{D/L},
\eqtn
and omit $S$ if $U=L$. When the explicit indication of~$S$ is not important, 
we use the term {\it generalized wreath product}. The generalized $S$-wreath product is  
{\it proper}  if $L\ne e$ and $U\ne D$. When $U=L$, the generalized $S$-wreath product coincides with the ordinary
wreath product.\medskip

Let $D=D_1D_2$, where $D_1$ and $D_2$ are trivially intersecting subgroups of $D$.
If $\cA_1$ and $\cA_2$ are S-rings over the groups $D_1$ and $D_2$ respectively, then the 
module 
$$
\cA=\Span\{\und{X_1\cdot X_2}: X_1\in\cS(\cA_1),\ X_2\in\cS(\cA_2)\}
$$
is an S-ring over the group $D$  whenever  $\cA_1$ and $\cA_2$ are commute with each other.
In this case, $\cA$ is called the {\it dot product} of $\cA_1$ and
$\cA_2$, and denoted by $\cA_1\cdot\cA_2$ \cite{LM98}. When $D=D_1\times D_2$, the
dot product coincides with the {\it tensor product} $\cA_1\otimes \cA_2$. The following 
statement was proved in \cite{EKP14}.

\lmml{050813b}
Let $\cA$ be an S-ring over an abelian group $D=D_1\times D_2$.  Suppose that $D_1$ and $D_2$
are $\cA$-groups. Then $\cA=\cA_{D_1}\otimes \cA_{D_2}$ whenever $\cA_{D_1}=\mZ D_1$ 
or $\cA_{D_2}=\mZ D_2$.
\elmm

The following two important theorems go back to Schur and Wielandt (see \cite[Ch.~IV]{Wie64}). 
The first of them is known as the Schur theorem on multipliers, see~\cite{EvdP09}. 

\thrml{261009b}
Let $\cA$ be an S-ring over an abelian group $D$. Then given an integer~$m$ coprime
to $|D|$, the mapping $X\mapsto X^{(m)}$, $X\in\cS(\cA)$, where
\qtnl{141014a}
X^{(m)}=\{x^m:\ x\in X\},
\eqtn
is a bijection. Moreover, $x\mapsto x^m$, $x\in D$, is a Cayley automorphism of~$\cA$.
\ethrm

Given a subset $X$ of an abelian group $D$, denote by $\tr(X)$ the {\it trace} of $X$, i.e. the union 
of all~$X^{(m)}$ over the integers $m$ coprime to~$|D|$. We say that 
$X$ is {\it rational} if $X=\tr(X)$. When $\tr(X)=\tr(Y)$ for some $Y\subseteq D$, the sets
$X$ and $Y$ are called {\it rationally conjugate}. For an S-ring $\cA$ over $D$, the module 
$$
\tr(\cA)=\Span\{\tr(X):\ X\in\cS(\cA\}
$$ 
is also an S-ring; it is called the {\it rational closure} of~$\cA$. Finally, the S-ring  is 
{\it rational} if it coincides with its rational closure, or equivalently, if each of its 
basic sets is rational.\medskip

In general, Theorem~\ref{261009b} is not true when $m$ is not coprime to the order of $D$. However, the following 
weaker statement holds.

\thrml{261009w}
Let $\cA$ be an S-ring over an abelian group $D$. Then given a prime divisor~$p$ of~$|D|$, 
the mapping $X\mapsto X^{[p]}$, $X\in 2^ D$, where
 \qtnl{030713a}
X^{[p]}=\{x^p:\ x\in X,\ \nmmod{|X\cap Hx|}{0}{p}\}
\eqtn
with $H=\{g\in D:\ g^p=1\}$, takes an $\cA$-set to an $\cA$-set.
\ethrm 

We complete the section by the theorem on separating subgroup that was proved 
in~\cite{EP05}.

\thrml{t100703}
Let $\cA$ be an   S-ring over a group $D$. Suppose that $X\in\cS(\cA)$ and $H\le D$ are such that 
$$
X\cap H\ne\emptyset\qaq X\setminus H\ne\emptyset\qaq
\grp{X\cap H}\le\rad(X\setminus H).
$$
Then $X=\grp{X}\setminus\rad(X)$ and $\rad(X)\le H\le\grp{X}$.
\ethrm

\section{S-rings and Cayley schemes}\label{111014a}
In this section, we freely use the language of association scheme theory; in our exposition,
we follow \cite{EP09,MP09}.

\sbsnt{The 1-1 correspondence.}
For a group $D$, denote by $R(D)$ the set of all binary relations on $D$ that are
invariant with respect to the group $D_{right}$ (consisting of the permutations of the set~$D$ induced
by the right multiplications in the group~$D$). Then the mapping
\qtnl{171014a}
2^D\to R(D),\quad X\mapsto R_D(X)
\eqtn
where $R_D(X)=\{(g,xg):\ g\in D,x\in X\}$, is a bijection. If $\cA$ is an S-ring
over the group $D$, then the pair 
\qtnl{210215a}
\cX=(D,S),
\eqtn 
where $S=R_D(\cS(\cA))$, is
an association scheme. Moreover, it is a {\it Cayley scheme} over $D$, i.e.
$D_{right}\le\aut(\cX)$. Each basis relation $s\in S$ of this scheme, is a Cayley digraph 
over~$D$ the connection set of which is equal to $es=\{g\in D:\ (e,g)\in s\}$. Conversely, 
given a Cayley scheme~\eqref{210215a} the module 
$$
\cA=\Span\{\und{es}:\ s\in S\}
$$ 
is an S-ring over~$D$.

\thrml{100909a}{\rm \cite{K85}}
The mappings $\cA\mapsto\cX$, $\cX\mapsto\cA$ form a 1-1 correspondence
between the S-rings and Cayley schemes over the group~$D$. 
\ethrm

It should be mentioned that the above correspondence preserves the inclusion. Moreover,
the mapping~\eqref{171014a} induces a ring isomorphism from $\cA$ onto the adjacency algebra
of the Cayley scheme $\cX$ associated with $\cA$. It follows that 
$c_{XY}^Z=c_{rs}^t$ for all $X,Y,Z\in\cS(\cA)$, where $r=R_D(X)$, $s=R_D(Y)$
and $t=R_D(Z)$. In particular, the number $|X|$ is equal to the valency $n_r$ of the relation $r$, and 
the S-ring~$\cA$ is commutative if and only if so is the Cayley scheme~$\cX$.

\sbsnt{Isomorphisms and schurity.}
We say that S-rings $\cA$ and $\cA'$ are (combinatorial) {\it isomorphic} if
the Cayley schemes associated with $\cA$ and $\cA'$ are isomorphic. Any isomorphism
between these schemes is called the {\it isomorphism} of $\cA$ and $\cA'$. The group 
$\iso(\cA)$ of all isomorphisms from $\cA$ to itself has a normal subgroup 
$$
\aut(\cA)=\{f\in\iso(\cA):\ R_D(X)^f=R_D(X)\ \text{for all}\; X\in\cS(\cA)\};
$$
any such $f$ is called a (combinatorial) {\it automorphism} of the S-ring~$\cA$.
In particular, if $\cA=\mZ D$ (resp. $\rk(\cA)=2$), then $\aut(\cA)=D_{right}$
(resp. $\aut(\cA)=\sym(D)$).\medskip

The S-ring $\cA$ is called {\it schurian} (resp. {\it normal}) if so is the Cayley scheme 
associated with $\cA$. Thus, $\cA$ is schurian if and only if $\cS(\cA)=\orb(\aut(\cA)_e,D)$, 
and normal if and only if $D_{right}\trianglelefteq\aut(\cA)$.\medskip

From our definitions, it follows that $\cA=\cA_1\wr\cA_2$ if and only if $\cX=\cX_1\wr\cX_2$
where $\cX$, $\cX_1$ and $\cX_2$ are the Cayley schemes associated with the S-rings~$\cA$, $\cA_1$ and
$\cA_2$ respectively. Similarly, $\cA=\cA_1\cdot\cA_2$ if and only if $\cX=\cX_1\otimes\cX_2$.
On the other hand, the tensor and wreath product of association schemes (and permutation groups) are special cases  
of the crested product introduced and studied in~\cite{BC05}. Thus, Theorem~\ref{211014a} below immediately follows from
Remark~23 of that paper and Theorems~21 and~22 proved there.

\thrml{211014a}
Let $\cA=\cA_1\ast\cA_2$, where $\ast\in\{\wr,\cdot\}$. Then  $\cA$ is schurian if and
only if so are $\cA_1$ and $\cA_2$. Moreover, 
$$
\aut(\cA_1\wr\cA_2)=\aut(\cA_1)\wr\aut(\cA_2)\qaq
\aut(\cA_1\cdot\cA_2)=\aut(\cA_1)\times\aut(\cA_2).
$$
\ethrm

The following simple statement is an obvious consequence of the definition of wreath product
and Theorem~\ref{211014a}.

\crllrl{021014a}
Let $\cA$ be an S-ring over a group $D$ and $H$ an $\cA$-group such that 
$\rk(\cA)=\rk(\cA_H)+1$. Then $\cA$ is isomorphic to the wreath product of $\cA_H$
by an S-ring of rank $2$ over the group $\mZ_{[D:H]}$.
Moreover, $\cA$ is schurian if and only if so is $\cA_H$.
\ecrllr

\sbsnt{Quasi-thin S-rings.} An S-ring $\cA$ is called {\it quasi-thin} if any of its basic sets consists of at most
two elements. Thus, $\cA$ is quasi-thin if and only if 
the Cayley scheme associated with $\cA$ is quasi-thin (the latter means that the valency 
of any its basic relation is at most~$2$).

\lmml{100914a}
Let $\cA$ be an S-ring over an abelian group $D$. Suppose that $X\in\cS(\cA)$ is such
that $|X|=2$ and $\grp{X}=D$. Then $\cA$ is quasi-thin.
\elmm
\proof The Cayley scheme~$\cX$ associated with $\cA$ is commutative,
because the group $D$ is abelian. Moreover, the relation $r=R_D(X)$ corresponding
to the set $X$, is of valency $|X|=2$. Thus, the equality $\grp{X}=D$ implies that $\cX$ is a $2$-cyclic
scheme generated by the tightly attached relation $r$ in the sense of~\cite{HM}. So, by Proposition~3.11 of that paper,
$\cX$ is a quasi-thin scheme. Therefore, the S-ring~$\cA$ is also quasi-thin.\bull

Following the theory of quasi-thin schemes in \cite{MP}, we say that a basic set $X\ne\{e\}$ 
of a quasi-thin S-ring~$\cA$ is an {\it orthogonal} if $X\subseteq Y\, Y^{-1}$ for some 
$Y\in\cS(\cA)$.

\lmml{100914y}
Any commutative quasi-thin  S-ring $\cA$ is schurian. Moreover, if
it has at least two orthogonals, then the group $\aut(\cA)_e$ has a faithful regular orbit.
\elmm
\proof The first statement immediately follows from ~\cite[Theorem~1.2]{MP}. To prove
the second one, denote by $\cX$ the Cayley scheme associated with the S-ring~$\cA$. Then
from \cite[Corollary~6.4]{MP} it follows that  the group $\aut(\cX)_{e,x}$ is trivial for
some $x\in D$. This means that $x^{\aut(\cX)_e}$ is a faithful regular orbit of the group $\aut(\cA)_e$.\bull

\section{S-rings over cyclic and dihedral groups}\label{231014v}

\sbsnt{Cyclic groups.}\label{141114b}
Let $C$ be a cyclic group of order $2^n$, $n\ge 1$. Then the group $\aut(C)$
consists of permutations $\sigma_m:x\mapsto x^m$, $x\in C$, where $m$ is an odd integer. In what
follows, $c_1$ denotes the unique involution in $C$.

\lmml{211113a}
Let $X\in\orb(K,C)$, where $K\le\aut(C)$. Then 
\nmrt
\tm{1} $\rad(X)=e$ if and only if $X$ is a singleton, or $n\ge 3$ and $X=\{x,\varepsilon x^{-1}\}$
where $x\in X$ and $\varepsilon\in \{e,c_1\}$,
\tm{2} if $K\ge\{\sigma_m:\ \mmod{m}{1}{2^{n-k}}\}$, then $2^k$ divides $|\rad(X)|$.
\enmrt
\elmm
\proof Statement (1) follows from \cite[Lemma~5.1]{EP02}, whereas statement~(2) is straightforward.\bull

Let $\cA$ be an S-ring over the group $C$. By the Schur theorem on multipliers,
the group $\rad(X)$ does not depend on a set $X\in\cS(\cA)$ that contains a generator of $C$.
This group is called the {\it radical} of $\cA$ and denoted by $\rad(\cA)$. Since $C$ is a $2$-group, from~\cite[Lemma~6.4]{EP02} 
it follows that  if $\rad(\cA)=e$, then either $n\ge 2$ and $\rk(\cA)=2$, or $\cA=\cyc(K,C)$, where $K\le\aut(C)$ is the group
generated by the automorphism taking a generator $x$ of $C$ to an element in $\{x,x^{-1},c_1x^{-1}\}$
(see also statement~(1) of Lemma~\ref{211113a}). In any case, $\cA$ is, obviously, schurian. In fact, the latter statement holds
for any S-ring over a cyclic $p$-group \cite{EKP}.\medskip 

For any basic set $X$ of the S-ring~$\cA$, one can form an $\cA$-section $S=\grp{X}/\rad(X)$. Then the radical
of the S-ring~$\cA_S$ is trivial. Since in our case, $|S|$ is a $2$-group, the result in previous paragraph shows
that either the S-ring $\cA_S$ is cyclotomic or $|S|$ is a composite number and $\rk(\cA_S)=2$. In the 
former case, $X$ is an orbit of an automorphism group of $C$, whereas in the  latter case, $X=\grp{X}\setminus\rad(X)$. Thus, any 
basic set of $\cA$ is either regular or equals the set difference of two distinct $\cA$-groups.

\lmml{021109a}
Let $\cA$ be a cyclotomic S-ring over a cyclic $2$-group.
Suppose that $\rad(\cA)=e$. Then $\rad(\cA_S)=e$ for any $\cA$-section $S$ such that $|S|\ne 4$.
\elmm
\proof Follows from \cite[Theorem~7.3]{EKP}.\bull

The following auxiliary lemma will be used in Section~\ref{140914b}.

\lmml{150414a}
Let $C$ be a cyclic $2$-group, and let $X$ and $Y$ be orbits of some subgroups of $\aut(C)$. Suppose
 that  $\grp{X}\ne\grp{Y}$  and $\rad(X)=\rad(Y)=e$. Then the product $X\,Y$ contains
no coset of $\grp{c_1}$.
\elmm
\proof Without loss of generality, we can assume that $\grp{Y}$ is a proper subgroup of~$\grp{X}$.
Then from statement (1) of Lemma~\ref{211113a}, it follows that
$X=\{x\}$ or $\{x,\varepsilon x^{-1}\}$, and $Y=\{y\}$ or $\{y,\varepsilon' y^{-1}\}$
where $|x|>|y|$. Therefore, the required statement trivially holds whenever $X$ or $Y$ is a singleton. Thus,
we can assume that $|X|=|Y|=2$, and hence  $|x|>|y|\ge 8$. Furthermore,
$$
X\,Y=\{xy,\quad \varepsilon'xy^{-1},\quad \varepsilon x^{-1}y,\quad \varepsilon''x^{-1}y^{-1}\},
$$
where $\varepsilon''=\varepsilon\varepsilon'$.  Suppose on the contrary, that this product contains a coset 
of $\grp{c_1}$. Then, obviously, $c_1xy=\varepsilon'xy^{-1}$ or  
$c_1\varepsilon x^{-1}y=\varepsilon''x^{-1}y^{-1}$.
In any case, $y^2\in\{e,c_1\}$ and hence $|y|\le 4$. Contradiction.\bull

\sbsnt{Dihedral groups.}\label{011014b}
Throughout this subsection, $D$ is a dihedral group and $C$ is the cyclic subgroup of $D$ such that all the 
elements in $D\setminus C$ are involutions. A set $X\subseteq D$ is called {\it mixed} if
the sets $X_0=X\cap C$ and $X\setminus X_0$ are not empty. For an
element $s\in D\setminus C$, we denote by $X_1=X_{1,s}$ the subset of 
$C$ for which $X=X_0\,\cup\,X_1s$. The following statement (in the other notation) was proved in~\cite{W49}.

\lmml{240913a}
Let $\cA$ be an S-ring over the dihedral group $D$ and $X$ a mixed basic set of $\cA$. Then
\nmrt
\tm{1} the sets $X_0$, $X_1s$ and $X$ are symmetric, and $X_0$ commutes with $X_1s$,
\tm{2} given an integer $m$ coprime to $|D|$, there exists a unique $Y\in\cS(\cA)$ such that $(X_0)^{(m)}=Y_0$.
\enmrt
\elmm

When it does not lead to confusion, the set $Y$ from statement~(2) of Lemma~\ref{240913a} will be also denoted by $X^{(m)}$; 
for $X\subseteq C$, this notation is consistent with~\eqref{141014a}.  For any $\cA$-set $X$ such that $X_0\ne\emptyset$, one can define 
its trace $\tr(X)$ to be the union of sets $X^{(m)}$, where $m$ runs over all integers coprime to $|D|$. The following statement was 
also proved in~\cite{W49}.

\lmml{290514b}
Let $\cA$ be an S-ring over the dihedral group $D$. Suppose that $X_0\ne\emptyset$ for all $X\in\cS(\cA)$. Then 
given an integer $m$ coprime to $|D|$, the mapping $X\mapsto X^{(m)}$ induces an algebraic isomorphism of~$\cA$; in particular,
$|X^{(m)}|=|X|$.
\elmm

The algebraic fusion of the S-ring $\cA$ from Lemma~\ref{290514b} with respect to the group of all algebraic  isomorphisms
defined in this lemma, is an S-ring any basic set of which is of the form $\tr(X)$ where $X\in\cS(\cA)$. This S-ring is called the 
{\it rational closure} of $\cA$ and denoted by $\tr(\cA)$. It should be stressed that this notation has sense only
if the hypothesis of Lemma~\ref{290514b} is satisfied.

\section{S-rings over $D=\mZ_2\times\mZ_{2^n}$: basic sets containing involutions}\label{231014a1}

In what follows, $C\le D$ is a cyclic group of order $2^n$ and $E$ is the Klein subgroup of~$D$. The non-identity elements of 
this subgroup are the involution $c_1\in C$ and the other two involutions $s\in D\setminus C$ and $c_1s$.
For an S-ring over~$D$ and an element $t\in E$, we denote by $X_t$ the basic set  that contains~$t$. Then
by the Schur theorem on multipliers (Theorem~\ref{261009b}),  this set is rational. In this section we completely 
describe the sets~$X_t$'s.

\thrml{250814a}
Let $\cA$ be an S-ring over the group $D$. Then the set $H=\bigcup_{t\in E}X_t$
is an $\cA$-group and for a suitable choice of $s$ one of the following statements holds
with $U=\grp{X_{c_1}}$:
\nmrt
\tm{1} $X_{c_1}=X_s=X_{sc_1}=U\setminus e$,
\tm{2} $X_{c_1}=U\setminus e$ and $X_s=X_{sc_1}=H\setminus U$,
\tm{3} $X_{c_1}=X_{sc_1}=U\setminus\grp{s}$ and $X_s=\{s\}$,
\tm{4} $X_{c_1}=U\setminus e$, $X_s=\{s\}$ and $X_{sc_1}=sX_{c_1}$.
\enmrt
\ethrm

The proof of Theorem~\ref{250814a}  will be given later. The following
auxiliary statement is, in fact, a consequence of the Schur theorem on multipliers. Below, we fix an S-ring~$\cA$ over
the group~$D$.

\lmml{290714c}
Suppose that $X\in\cS(\cA)$ contains two elements $x$ and $y$ such that $|x| > |y|\geq 2$. Then $x\{e,c_1\}\subseteq X$.
\elmm
\proof Set $m=1 + |x|/2$. Then by Schur's theorem  on multipliers, $Y:=X^{(m)}$ is a basic set of $\cA$. On the other hand,
since $|x| > |y|$, we have $y^m=y$. Thus, $y^m\in X$. This implies that $X=Y$, and hence $x^m\in X$. Since $|x|>2$ and 
$x^m=xc_1$, we conclude that $x\{e,c_1\}\subseteq X$ as required.\bull\medskip

In the following lemma, we keep the notation of Theorem~\ref{250814a}.

\lmml{290714b} 
Either $X_{c_1}=U\setminus e$, or $H$ is an $\cA$-group and statement~(3) of Theorem~\ref{250814a} holds. 
\elmm
\proof
The statement is trivial if the set $X:=X_{c_1}$ is contained in~$E$. So, we can assume that $X$ contains at least one element 
of order greater than two. Then $x\{e,c_1\}\subseteq X$ for each $x\in X$ with $|x|>2$ (Lemma~\ref{290714c}).
Thus, $c_1\in\rad(X\setminus E)$. We observe that $X\cap E$ is equal to one of the following sets:
$$
\{c_1\}\qoq \{c_1,s,sc_1\}\qoq \{c_1,t\},
$$
where $t\in\{s,sc_1\}$.
In the first two cases, $c_1\in\rad(X\setminus\{c_1\})$. So by Theorem~\ref{t100703} with $H=\grp{c_1}$, we 
conclude that $X=\grp{X}\setminus\rad(X)$ and $c_1\not\in\rad(X)$. However, the only non-trivial subgroups of $D$ 
not containing $c_1$, are $\grp{s}$ and $\grp{sc_1}$. Since none of them equals $\rad(X)$,
we conclude that $\rad(X)=e$ and $X=U\setminus e$.\medskip

In the remaining case, $X\cap E=\{c_1,t\}$, and hence $|c_1 X\cap X|=|X|-2$. Since $c_1\in X$, this implies that  the latter number 
equals~$c_{XX}^X$ (see Section~\ref{120813a}). 
Therefore,  $|x^{-1}X\cap X|=|X|-2$ for each $x\in X$. On the other hand, the set $\{c_1,t\}t=\{c_1t,e\}$ does not
intersect $X$. Thus, $t(X\setminus E)=X\setminus E$, and hence 
$$
\rad(X\setminus E)=E.
$$
By Theorem~\ref{t100703} with 
$H=E$, we conclude that $X=\grp{X}\setminus\rad(X)$ and $E\setminus\rad(X)$ is contained in~$X$. Therefore, $\rad(X)=\grp{t'}$,
where $t'$ is the element of $\{s,sc_1\}$ other than $t$. Thus, $X=U\setminus\grp{t'}$, $H=U$ is an $\cA$-group and we are done.\bull\medskip

{\bf Proof of Theorem~\ref{250814a}.} By Lemma~\ref{290714b}, we can assume that $X:=X_{c_1}=U\setminus e$. 
If, in addition, $X$ contains $s$ or $sc_1$, then $H=U$ is an $\cA$-group and statement~(1) of the theorem holds. Thus, we can 
also assume that 
$X\ne X_t$ for each $t\in \{s,sc_1\}$. Suppose first that $X_s=X_{sc_1}$; denote this set by $Y$. Then $Y\cap E=\{s,sc_1\}$, and hence
$c_1\in\rad(Y)$ (Lemma~\ref{290714c}). Since $X=U\setminus e$, this implies that $U\le\rad(Y)$. 
Thus, $H=\grp{Y}$ is an $\cA$-group, $X_s=X_{sc_1}=Y=H\setminus U$ and statement~(2) holds.\medskip

Let now $X_s\ne X_{sc_1}$ and $t\in \{s,sc_1\}$. Then $Y\cap E=\{t\}$ where $Y=X_t$. It follows that $|c_1 Y\cap Y|=|Y|-1$, and hence
\qtnl{280814a}
\und{Y}^2 = |Y|e+(|Y|-1) \und{X} + \cdots,
\eqtn
where the omitted terms in the right-hand side contain neither $e$ nor elements of $X$ with non-zero coefficients. However,
$|X|c_{YY}^X=|Y|c_{YX}^Y$ because $X=X^{-1}$ and $Y=Y^{-1}$.  Since $c_{YY}^X=|Y|-1$, this implies
that $|X|$ is divided by $|Y|$. If, in addition, $|Y|=1$, then we have $\{X_s,X_{sc_1}\}=\{\{t\},tX_{c_1}\}$ and $H=U\cup tU$ is an $\cA$-group. So statement~(4) holds.   If $|Y|\ne 1$, equality~\eqref{280814a}
implies that $|Y|=|X|=|U|-1$ and $\und{Y}^2 = |Y|e+(|Y|-1) \und{X}$. It follows that $c_{YX}^Y=|X|-1$. Therefore,
$$
|tU\cap Y|=|X|-1=|U|-2.
$$
This is true for $t=c_1$ and $t=sc_1$. On the other hand, $sU=sc_1U$ and the sets $X_s$ and $X_{sc_1}$ are disjoint. Thus, 
$$
|U|=|tU|\ge |sU\cap X_s|+|sc_1U\cap X_{sc_1}|=2(|U|-2).
$$
It follows that $|U|=2$ or $|U|=4$. In the former case, $H=E$ and statement~(4) trivially holds. In the
latter one, $|X|=|X_s|=3$ and $|sU\cap Y|=2$. Therefore, there exists a unique $x\in X_s$ outside $sU$. It
follows that $|xU\cap X_s|=1$ which is impossible by Lemma~\ref{090608a}.\bull\medskip

\section{S-rings over $D=\mZ_2\times\mZ_{2^n}$: non-regular case}\label{140914a}

A set $X\subseteq D$ is said to be {\it highest} (in $D$) if it contains an element of order $2^n$.
Given an S-ring $\cA$  over $D$, denote by $\rad(\cA)$ the group generated by
the groups $\rad(X)$, where $X$ runs over the highest basic sets of $\cA$.
Clearly, $\rad(\cA)$ is an $\cA$-group, and it is equal to $e$ if and only if each highest basic set of $\cA$ has trivial radical. 
In what follows, we say that $\cA$ is {\it regular}, if  each highest basic set of $\cA$ is regular. Now, the 
main result of this section can be formulated as follows.

\thrml{170414a}
Let $\cA$ be an S-ring over the group $D$. Suppose that $\rad(\cA)=e$. Then $\cA$ is either regular
or rational. Moreover, in the latter case, $\cA=\cA_H\otimes\cA_L$, where $\rk(\cA_H)=2$ and  
$|L|\le 2\le |H|$; in particular, $\cA$ is schurian.
\ethrm

The proof of Theorem~\ref{170414a} will be given in the end of the section. The key point of the proof is the
following statement.

\thrml{290714a}
Let $\cA$ be  an S-ring over the group $D$. Then any non-regular basic set of $\cA$ either intersects $E$,
or has a non-trivial radical.
\ethrm
\proof
Let $X$ be a non-regular basic set of $\cA$ that does not intersect $E$. Then the minimal order of an element in $X$
equals $2^m$ for some $m\ge 2$. Denote by $X_m$ the set of all elements in~$X$ of order $2^m$. Clearly,
each of the sets  $X\setminus X_m$ and $X_m$ is non-empty. Suppose, towards a contradiction, that $\rad(X)=e$. Then $c_1\not\in\rad(X)$, and $c_1\in\rad(X\setminus X_m)$ (Lemma~\ref{290714c}). It follows that
\qtnl{300814u}
c_1\not\in\cA,
\eqtn
because otherwise $c_1X$ is a basic set other than $X$ that intersects $X$.\medskip

Denote by $K$ the setwise stabilizer of $X$ in the group $G\cong\mZ_{2^n}^*$ of all permutations $x\mapsto x^m$, 
$x\in D$, where $m$ is an odd integer.  Then
by Schur's theorem on multipliers, $X_m$ is the union of at most  two $K$-orbits (one inside $C$ and the other outside). The 
radicals of these orbits must be trivial, because
$\rad(X)=e$ and $c_1\in\rad(X\setminus X_m)$. Thus, by statement~(1) of Lemma~\ref{211113a}, we have 
\qtnl{300814a}
X_m=\{x\}\qoq \{x,x^{-1}\varepsilon\}\qoq\{x,ys\}\qoq\{x,x^{-1}\varepsilon,ys,y^{-1}\varepsilon s\},
\eqtn
where $x,y\in X_m$ are such that $\grp{x}=\grp{y}$, and $\varepsilon\in\{e,c_1\}$. It should be mentioned that $x\ne \varepsilon y$, 
for otherwise $\varepsilon s\in\rad(X)$.\medskip

Let us define $\cA$-groups $U$ and $H$ as in Theorem~\ref{250814a}. Then 
$X\subseteq D\setminus H$, because $X$ does not intersect $E$. By the definition of $H$, this implies that it does not contain elements
of order $2^m$. Since $U\le H$, we conclude that 
$xU\cap X\subseteq X_m$ for each $x\in X_m$. Therefore, $X_m$ is a disjoint union
of some sets $xU\cap X$ with such $x$. However, by Lemma~\ref{090608a} the number $\lambda:=|xU\cap X|$ doesn't depend 
on a choice of $x\in X$. Thus, $\lambda$ divides $|X_m|$.
By~\eqref{300814a} this implies that
$\lambda\in\{1,2,4\}$. Moreover, setting $Y$ to be the basic set 
containing $c_1$, we have
\qtnl{141114a}
c_{X^{} Y^{}}^X=\lambda-1,
\eqtn
because $U=Y\setminus e$ or $U=Y\setminus \grp{\varepsilon s}$, 
and $x\ne \varepsilon y$.\medskip

Denote by $\alpha$ the number of $z\in X_m$ for which $c_1z\not\in X_m$.
If $\alpha=1$, then from Theorem~\ref{261009w}, it follows that 
$\cS(\cA)$ contains $X^{[2]}=\{z^2\}$ for an appropriate $z\in X_m$. 
Since $z\not\in E$, this implies that $c_1\in\cA$ in contrast
to~\eqref{300814u}. Thus, $\alpha\ne 1$. Therefore,  $\alpha$ is an 
even number less or equal than $|X_m|\le 4$ (see~\eqref{300814a}).
Moreover, it is not zero, because otherwise $c_1\in\rad(X)$. Besides,
from~\eqref{141114a} it follows that
\qtnl{310714a}
|X|\,(\lambda - 1)=|X|c_{X^{} Y^{}}^X = 
|Y|\,c_{X^{} X^{-1}}^Y = |Y|\,|c_1X\cap X| = |Y|\,(|X|-\alpha).
\eqtn
Since $|X|>|X_m|\ge\alpha$, this implies that the right-hand side of 
the equality is not zero. Thus, $\lambda\ne 1$, and finally
\qtnl{010914a}
\lambda,\alpha\in\{2,4\}.
\eqtn
\lmml{310714b} 
In the above notation $|X|\geq 2|X_m|$, and the equality holds only if  
\nmrt
\tm{1} $X_m$ is a union of two $K$-orbits and $X\setminus X_m$ is a $K$-orbit,
\tm{2} any element in $X\setminus X_m$ is of order $2^{m+1}$.
\enmrt
\elmm
\proof By the Schur theorem on multipliers, the stabilizers of 
an element $x\in X$ in the groups $K$ and $G$, coincide. However, 
the stabilizer in $G$ consists of raising to power $1+i|x|$, where
$i=0,1,\ldots,2^n/|x|-1$. Therefore, $|K_x|=2^n/|x|$. For $x\in X_m$ and 
$y\in X\setminus X_m$, this implies that
$|K_x|\ge 2|K_y|$, and hence
$$
|x^K|=\frac{|G|}{|K_x|}\le \frac{|G|}{2|K_y|}=\frac{|y^K|}{2}.
$$
Taking into account that $X_m$ is a disjoint union of 
at most two $K$-orbits, we obtain that
$$
|X|-|X_m|\geq |y^K|\geq 2|x^K|\geq |X_m|
$$  
as required. Since the equality holds only if the second and third
inequalities in the above formula are equalities, we are done.\bull\medskip

We observe that $|Y|\ne\lambda-1$: indeed, for $\lambda=2$, this follows from~\eqref{300814u} whereas
for $\lambda=4$, the assumption $|Y|=\lambda-1$ implies by \eqref{310714a} an impossible equality $|X|=|X|-\alpha$.
Thus, by~\eqref{310714a} and Lemma~\ref{310714b} 
we have
$$
\frac{\alpha\,|Y|}{|Y|-(\lambda-1)}=|X|\ge 2|X_m|\ge 2\alpha.
$$
Furthermore, if $\lambda=2$, then $|Y|=2$ and $|X|=2\alpha$. On the other hand, if $\lambda=4$, then $\lambda-1<|Y|\le 6$ and
$|Y|\in \{2^a-1,2^a-2\}$ for some~$a$ (Theorem~\ref {250814a}); but then $|Y|=6$ and $|T|=2\alpha$. 
Thus, by~\eqref{010914a} there are exactly four possibilities:
\nmrt
\tm{1} $\alpha =2, \lambda=2, |Y|=2, |X|=4, |X_m|=2$,
\tm{2} $\alpha =4, \lambda=2, |Y|=2, |X|=8, |X_m|=4$,
\tm{3} $\alpha =2, \lambda=4, |Y|=6, |X|=4, |X_m|=2$,
\tm{4} $\alpha =4, \lambda=4, |Y|=6, |X|=8, |X_m|=4$.
\enmrt
In all cases, $|X|=2|X_m|$. Therefore, by Lemma~\ref{310714b} and \eqref{300814a}, we conclude that 
$X_m=\{x,ys\}$ in cases (1) and~(3), and  $X_m = \{x,x^{-1}\varepsilon,sy,sy^{-1}\varepsilon\}$ 
in cases (2) and~(4). Moreover, since the number $|Y|$ is even, $U\setminus Y$ is a group of order two. Without loss 
of generality, we assume that it is $\grp{s}$.\medskip

Let $\pi$ be the quotient epimorphism from $D$ to $D'=D/U$. Then the group $D'$ is cyclic, the S-ring $\cA'=\cA_{D'}$ is 
circulant\footnote{Any S-ring over a cyclic group is called a circulant one.} and $X'=\pi(X)$ is a non-regular basic set of it (the elements in $X'_m=\pi(X_m)$ and in $X_{}'\setminus X'_m$ have
different orders). However, any non-regular basic set of a circulant S-ring over a $2$-group is a set difference of two its subgroups
(see Subsection~\ref{141114b}). Therefore,
$$
X'=\grp{X'}\setminus\rad(X').
$$
Since $X_{}'\ne X'_m$, this implies that $|X'|\ge 3$. On the other hand, $|X'|=|X|/\lambda$ by the definition of $\lambda$.  Thus,
we can exclude cases (1), (3) and (4). In case (2) let $|\rad(X')|=2^i$ for some $i\ge 0$. Then $4=|X'|=2^{i+2}-2^i=3\,2^i$.
Contradiction.\bull\medskip

Any basic set $X$ of an S-ring $\cA$ over $D$ that intersects the 
group~$E$, must contain an involution. Therefore, such $X$ is rational. 
By Theorem~\ref{290714a}, this proves the following statement.

\crllrl{011213d}
Let $X$ be a basic set of an S-ring  over $D$. Suppose that $\rad(X)=e$. Then $X$
is either regular or rational.\bull
\ecrllr

{\bf Proof of Theorem~\ref{170414a}.}\ Suppose that $\cA$ is not regular. Then there exists a highest set $X\in\cS(\cA)$ 
that is not regular.
Since $\rad(X)=e$, we conclude by Theorem~\ref{290714a} that the set $X\cap E$ is not empty. Therefore,
$X$ is contained in the $\cA$-group~$H\ge E$ defined in Theorem~\ref{250814a}. But then $H=D$, because
the set $X$ is highest. Now, the first statement follows, because the
S-ring $\cA_H=\cA$ is, obviously, rational. Moreover, statements (2) and (3) of Theorem~\ref{250814a} do not
hold, because $\rad(\cA_H)=\rad(\cA)=e$. Thus, the second statement of our theorem is true for $L=e$ 
(resp. $L=\grp{s}$) if statement~(1) (resp. statement~(4)) of Theorem~\ref{250814a} holds.\bull\medskip

From the proof of Theorem~\ref{170414a}, it follows that if one of the highest basic sets of $\cA$ is not regular,
then all highest basic sets are rational. This implies the following statement.

\crllrl{180514a}
Let $\cA$ be an S-ring over the group $D$. Suppose that $\rad(\cA)=e$. Then either every highest basic set of
$\cA$ is regular, or every highest basic set of $\cA$ is rational.\bull
\ecrllr

\section{S-rings over $D=\mZ_2\times\mZ_{2^n}$: regular case}\label{100414a}

Throughout this section, $C=C_n$ is a cyclic subgroup
of~$D=D_n$ that is isomorphic to $\mZ_{2^n}$.
We denote by $c_1$, $c_2$ and $s$, respectively, the unique involution in $C$, one of the two elements of $C$ of order~$4$, one of the two
involutions in $D\setminus C$.  The main result is given by the following theorem.

\thrml{230214a}
Let $\cA$ be a regular S-ring  over the group $D$. Suppose that $\rad(\cA)=e$. Then $\cA$ is a cyclotomic S-ring. 
More precisely, $\cA=\cyc(K,D)$, where $K\le\aut(D)$ is one of the groups listed in Table~\ref{280414b}.
\ethrm
\begin{table}
	\centering
		\begin{tabular}[tc]{|c|l|c|c|c|}
\hline 
$K$ & generators & $|K|$ & $n$ & comment \\
\hline
$K_1$ & $(x,s)\mapsto (x,s)$ & 1 & $n\ge 2$ & $X_1=\emptyset$\\
\hline
$K_2$ & $(x,s)\mapsto (x^{-1},s)$ & 2&  $n\ge 3$ & $X_1=\emptyset$ \\
\hline
$K_3$ & $(x,s)\mapsto (c_1x^{-1},s)$ &  2& $n\ge 3$ & $X_1=\emptyset$ \\
\hline
$K_4$ & $(x,s)\mapsto (x^{-1},sc_1)$      & 2  &  $n\ge 3$  & $X_1=\emptyset$ \\
\hline
$K_5$ & $(x,s)\mapsto (c_1x^{-1},sc_1)$      & 2  &  $n\ge 3$  & $X_1=\emptyset$ \\
\hline
$K_6$ & $(x,s)\mapsto (sc_2x,sc_1)$,\ $(x,s)\mapsto (x^{-1},s)$   & 4 &      $n\ge 4$  & $X_a\ne\emptyset$\\
\hline
$K_7$ & $(x,s)\mapsto (sc_2x,sc_1)$,\ $(x,s)\mapsto (c_1x^{-1},s)$             & 4 & $n\ge 4$  & $X_a\ne\emptyset$\\
\hline
$K_8$ & $(x,s)\mapsto (sx^{-1},s)$ & 2 &$n\ge 4$  & $X_a\ne\emptyset$\\
\hline
$K_9$ & $(x,s)\mapsto (sc_1x^{-1},s)$ & 2&  $n\ge 4$  & $X_a\ne\emptyset$\\
\hline
$K_{10}$ & $(x,s)\mapsto (sc_2x,sc_1)$ &  2& $n\ge 3$  & $X_a\ne\emptyset$\\
\hline
$K_{11}$ & $(x,s)\mapsto (sc_2x^{-1},sc_1)$ &  2& $n\ge 4$  & $X_a\ne\emptyset$\\
\hline
		\end{tabular}\\[2mm]
	\caption{The groups of cyclotomic rings with trivial radical}
	\label{280414b}
\end{table}
\crllrl{100514d}
Under the assumptions of Theorem~\ref{230214a}, let $K=K_i$, where $i=1,\ldots,11$. Then the following 
statements hold:
\nmrt
\tm{1} $\grp{c_1}$ is an $\cA$-group,
\tm{2} $C$ is an $\cA$-group if and only if $i\le 5$,
\tm{3} $\grp{\varepsilon s}$ with $\varepsilon\in\{e,c_1\}$, is an $\cA$-group if and only if
$i\in\{1,2,3,8,9\}$.\bull
\enmrt
\ecrllr

In what follows, given a basic set $X\in\cS(\cA)$, we denote by $X_0$ and $X_1$ the uniquely determined subsets of $C$ for which
$X=X_0\cup sX_1$.\medskip

{\bf Proof of Theorem~\ref{230214a}.} Let $X$ be a highest basic set 
of the S-ring~$\cA$. Then it is regular by the theorem hypothesis.
By the Schur theorem on multipliers, this implies that 
if the set $X_a$ is not empty for some $a\in\{0,1\}$, then $X_a$ is an orbit of 
the group $\aut(C)_{\{X_a\}}$. Therefore, $X_a$ is of the form given in
statement~(1) of Lemma~\ref{211113a}. The rest of the proof consists of Lemmas~\ref{040914a},
\ref{040914b} and \ref{040914c} below:  in the first  one, $X_0$ or $X_1$ is empty, and in the 
other two, both $X_0$ and $X_1$ are not empty and $|X_0|=2$ or $1$,
respectively.

\lmml{040914a}
Let $X_0=\emptyset$ or $X_1=\emptyset$. Then $\cA=\cyc(K_i,D)$ with $i\in\{1,2,3,4,5\}$.
\elmm
\proof
 Without loss of generality, we can assume that $n\ge 3$ and  $X_1=\emptyset$. Then $X=X_0$ generates $C$. Therefore, $C$ is 
an $\cA$-group and $X$ is a 
highest basic set of a circulant S-ring $\cA_C$. Since $\rad(X)=e$, this implies that $\rad(\cA_C)=e$. If, in addition, 
$\grp{s}$ is an $\cA$-group, then $\cA=\cA_C\otimes\cA_{\grp{s}}$
by Lemma~\ref{050813b}, and hence $\cA=\cyc(K_i,D)$ with $i=1,2,3$. Thus, we can assume  that
\qtnl{290414a}
s\not\in\cA.
\eqtn

Let us prove by induction on $n$ that $\cA=\cyc(K_i,D)$ with $i=4$ or $5$. For $n=3$ this statement can be 
verified by a computer computation. Let $n>3$. Denote by $X'$ the basic set of $\cA$ that contains $x'=xs$,
where $x\in X$ is a generator of~$C$. Then by the theorem hypothesis, 
$X'$ is a regular set with trivial radical. It follows that $C'=\grp{X'}$ is the order $2^n$ cyclic subgroup of $D$ 
other than $C$. In particular,
$$
\rad(\cA_{C'})=\rad(X')=e.
$$  

Besides, since $C^2=(C')^2$, the S-rings $\cA_C$ and $\cA_{C'}$ have the same basic sets inside the group $C^2$; in particular,
$|X|=|X'|$. Moreover, these S-rings are not Cayley isomorphic. Indeed, otherwise $X'=sY$, where $Y=X^{(m)}$ for some
odd~$m$. Then $s$ is the only element that appears in the product $\und{Y}^{-1}\und{X'}$ with multiplicity~$|X|$. However, 
in this case $s\in\cA$, contrary  to~\eqref{290414a}.  Thus,  
\qtnl{040914f}
X=\{x,\varepsilon x^{-1}\}\qaq X'=\{sx,sc_1\varepsilon x^{-1}\}.
\eqtn
Set $i=4$ or $5$ depending on $\varepsilon=e$ or $c_1$, respectively. Then 
from~\eqref{040914f} and the Schur theorem on multipliers, it follows that the S-rings $\cA$ and $\cyc(K_i,D)$ have the same 
highest basic sets. We also observe that $D_{n-1}$ is an $\cA$-group.\medskip

Since the S-ring $\cA_C$ is cyclotomic, it follows from~\eqref{040914f} that  $Y=\{x^2,x^{-2}\}$
is a basic set of $\cA$.  Denote by $Y'$ the basic set containing $sx^2$. Then, obviously, 
$$
Y'\subseteq X\,X'=\{sx^2,sc_1x^{-2},s,sc_1\}.
$$
However, $|sx^2|\ge 8$, because $n\ge 4$.  Thus, $Y'\subseteq\{sx^2,sc_1x^{-2}\} $: otherwise $Y\cap E$ is
not empty and from Theorem~\ref{250814a} it follows that $|Y'|>4$.  Therefore,
 $Y'$ is regular and $\rad(Y')=e$.  Since also $\rad(Y)=e$
and $Y,Y'$ are highest basic sets of the S-ring $\cA_{D_{n-1}}$, the latter 
satisfies the hypothesis of Lemma~\ref{040914a}. By the induction, we conclude that
$\cA_{D_{n-1}}=\cyc(K_4,D_{n-1})$. Thus, $\cA=\cyc(K_i,D)$ as required.\bull

\lmml{040914b}
Let $|X_0|=2$ and $X_1\ne\emptyset$. Then $\cA=\cyc(K_i,D)$ with $i\in\{6,7\}$.
\elmm
\proof By statement~(1) of Lemma~\ref{211113a}, we have $X_0=\{x,\varepsilon x^{-1}\}$.  Since $X$ is regular, this implies that
\qtnl{010514a}
X=\{x,\varepsilon x^{-1},sy,s\varepsilon y^{-1}\}
\eqtn
for some generator $y$ of the group $C$. For $n=3$, we tested in the computer that no S-ring over $D$  
has a highest basic set $X$ such that $X_0=\{x,\varepsilon x^{-1}\}$ and $X_1\ne\emptyset$.
Suppose that $n\ge 4$. Let us  prove that the lemma statement holds for $i=6$ or $i=7$ depending on 
$\varepsilon=e$ or $\varepsilon=c_1$, respectively. For $n=4$, we tested this statement in the computer. Thus, 
in what follows, we can assume that $n\ge 5$. 

\lmml{010514b}
In the above notations, the following statements hold:
\nmrt
\tm{1} $C_{n-1}$ is an $\cA$-group whereas $\grp{s}$ and $\grp{sc_1}$ are not,
\tm{2} $Y_x=\{x^{\pm 2},c_1x^{\pm 2}\}$ and $Z_x=\{s x^{\pm 2}\}$ are $\cA$-sets,
\tm{3} $y=xc_2$ for a suitable choice of $y$ and $c_2$.\footnote{One can interchange $y$ and $y^{-1}$,
and $c_2$ and $c_2^{-1}$.}
\enmrt
\elmm
\proof Since $n\ge 5$, we have $x^2\ne x^{-2}$ and $y^2\ne y^{-2}$. Besides, since neither $s$
nor $sc_1$ belongs to $\rad(X)$, we have $x^2\ne y^{\pm 2}\ne x^{-2}$. Thus,
\qtnl{030514a}
|\{x^2, x^{-2}, y^2, y^{-2}\}|=4.
\eqtn
However, $X^{[2]}=\{x^2, x^{-2}, y^2, y^{-2}\}$ and $X^{[2]}$ is an $\cA$-set by Theorem~\ref{261009w}.
Thus, the first part of statement~(1) holds, because $C_{n-1}=\grp{X^{[2]}}$. To  prove the second part of 
statement~(1), suppose on the contrary that $L:=\grp{s\varepsilon'}$ is an $\cA$-group
for some $\varepsilon'\in\{e,c_1\}$. Then the circulant S-ring $\cA_{D/L}$ has a basic set 
$$
\pi(X)=\{\pi(x),\pi(\varepsilon x^{-1}),\pi(y),\pi(\varepsilon y^{-1})\},
$$ 
where $\pi:D\to D/L$ is the quotient epimorphism. However, from~\eqref{030514a} it easily follows that
$|\rad(\pi(X))|\ge 2$. Therefore, one of the quotients $\pi(x)/\pi(\varepsilon x^{-1})$,  $\pi(x)/\pi(y)$ or  
$\pi(x)/\pi(\varepsilon y^{-1})$ has order~$2$. This implies, respectively, that the order of $\pi(x)$ is $8$, 
$\pi(x)=\pi(c_1)\pi(y)$, and $\pi(x)=\pi(c_1)\pi(\varepsilon y^{-1})$. The former case is impossible, because $n\ge 5$,
 whereas in the other two, we have $x\in c_1yL$ and $x\in c_1\varepsilon y^{-1}L$ which is impossible due 
to~\eqref{030514a}.\medskip

To prove the second part of statement~(2), we observe that  $C_{n-1}\cup\tr(X)$ is an $\cA$-set. But the complement to it  
in~$D$ coincides with $sC_{n-1}$. Thus, it is also an $\cA$-set. Besides,
\qtnl{110214a}
\und{X}^2\,\circ\,\und{sC_{n-1}}=2\,\und{sX'},
\eqtn
where $X'=\{(xy)^{\pm 1},\varepsilon(xy^{-1})^{\pm 1}\}$. Therefore, $sX'$ is an $\cA$-set.
However, it is easily seen that $|xy|\ne|\varepsilon xy^{-1}|$. Moreover, $|xy|=2^{n-1}$ or 
$|\varepsilon xy^{-1}|=2^{n-1}$, because 
$x$ and $y$ are generators of the group $C$ and $n\ge 3$. Thus, the
elements $sxy$ and $s\varepsilon xy^{-1}$ can not belong to the same basic set of $\cA$.
Indeed, otherwise assuming $|xy|=2^{n-1}$, we conclude by Lemma~\ref{290714c} that this basic set contains $sxyc_1$.
Then $xyc_1\in X'$, and hence $xyc_1\in\{(xy)^{-1},\varepsilon x^{-1}y\}$. Consequently, $(xy)^2=c_1$ or $x^2=c_1\varepsilon$. In any
case, $n-1\le 2$. Contradiction. A similar argument leads to a contradiction when $|\varepsilon xy^{-1}|=2^{n-1}$.
In the same way, one can verify that no two elements, one in $\{s(xy)^{\pm 1}\}$ and the other one in $\{s\varepsilon (xy^{-1})^{\pm 1}\}$,
cannot belong to the same basic set of $\cA$. Thus, $sX'$ is a disjoint union
of two $\cA$-sets of the form $\{sz^{\pm 1}\}$ with $z\in C$, and one of them consists of elements 
of order $2^{n-1}$. This implies that $Z_x$ is an $\cA$-set, as required.\medskip

 To prove statement~(3), suppose on the contrary that $x^4\ne y^{\pm 4}$. Then since $Y:=X^{[2]}$, $Z_x$, and $Z_y$ 
are $\cA$-sets, the S-ring $\cA$ contains the element
$$
(\und{Y_{}}\,\und{Z_x})\,\circ\, (\und{Y_{}}\,\und{Z_y})
=2s(2e+x^{\pm 2}y^{\pm 2}).
$$
Since $n\ge 3$, this implies that only $s$ appears in the right-hand 
side with multiplicity~$4$. By the Schur-Wielandt principle, this implies that $s\in\cA$. However, this contradicts the second part of statement~(1).\medskip

To complete the proof,  we note that by statement~(3), we have $Y_x=X^{[2]}$. Since $X^{[2]}$ is an $\cA$-set,
the first part of statement~(2) follows.\bull

Let us continue the proof of Lemma~\ref{040914b}.
Denote by $\cA_i$ the minimal S-ring over $D$ that contains $X$ as a basic set; we recall that
$i=6$ or $i=7$ depending on $\varepsilon=e$ or $\varepsilon=c_1$. We claim that
\qtnl{040514c}
\cA_i=\cyc(K_i,D).
\eqtn
Indeed, from statements~(2) and~(3) of Lemma~\ref{010514b}, it follows that the sets
$X$, $Y_x$, and $Z_x$ are orbits of the group $K_i$ (see Table~\ref{tblsg1}, where the generic orbits of the group $K_i$ 
contained inside $C_{n-1}$  and $sC_{n-1}$ are given).
\begin{table}[h]
\begin{center}
\begin{tabular}{|l|l|}
\hline 
     $C_{n-1}$   & $sC_{n-1}$       \\
\hline 
$x^{\pm 2},c_1x^{\pm 2}$ &  $s\varepsilon x^{\pm 2}$     \\
\hline 
  $x^4,x^{-4}$ &  $sx^{\pm 4},sc_1x^{\pm 4}$     \\
\hline 
  $\cdots$   &  $\cdots$                 \\
\hline 
$c_2^{},c_2^{-1}$ &  $sc_2^{},sc_2^{-1}$  \\
\hline 
  $c_1$   &  $s,sc_1$  \\	
\hline 
  $e  $   &  \\	
\hline 
\end{tabular}\\[2mm]
\caption{The orbits of $K_6$ and $K_7$}\label{tblsg1}
\end{center}
\end{table}
Therefore, by the Schur theorem on multipliers, we have  
\qtnl{040514a}
(\cA_i)_{D\setminus D_{n-2}}=\cyc(K_i,D)_{D\setminus D_{n-2}}.
\eqtn
Next, $C':=\grp{Z_x}$ is a cyclic $\cA$-group of order $2^{n-1}$ other than $C_{n-1}$. Moreover, $Z_x$ is a highest
basic set of the circulant S-ring $(\cA_i)_{C'}$. Therefore, this S-ring has trivial radical. From the results discussed
just after Lemma~\ref{211113a}, it follows that it is the cyclotomic S-ring $\cyc(K',C')$, where
$K'$ is the subgroup of $\aut(C')$ that has $Z_x$ as an orbit. Since $\orb(K_i,C')=\orb(K',C')$,
we conclude that
\qtnl{040514b}
(\cA_i)_{C_{n-2}}=\cyc(K_i,D)_{C_{n-2}}.
\eqtn
Let us complete the proof of~\eqref{040514c}. To do this, taking into account that $\varepsilon(e+c_1)=e+c_1$, we find that 
\qtnl{040514u}
\und{Y_x}\,\und{Z_x}=sx^{\pm 4}(e+c_1)+2s(e+c_1).
\eqtn
Moreover, since $n\ge 5$, the elements $x^{\pm 4},sx^{\pm 4}c_1$
appear in the right-hand side with coefficient~$1$. By the Schur-Wielandt principle, this implies that $\{s,sc_1\}x^{\pm 4}$ 
and $\{s,sc_1\}$ are $\cA$-sets. Thus,
\qtnl{040514e}
(\cA_i)_{sC_m\setminus sC_{m-1}}=\cyc(K_i,D)_{sC_m\setminus sC_{m-1}}\qaq
(\cA_i)_{sC_1}=\cyc(K_i,D)_{sC_1},
\eqtn
where $m=n-2$. For all $m=n-3,\ldots,2$ the first equality is proved in a similar way by the induction on~$m$; the sets $Y_x$ and 
$Z_x$ in equation~\eqref{040514u} are  replaced  by the $\cA$-sets $\{x^{\pm 2^{m+1}}\}$ and $(e+c_1)\{x^{\pm 2^{m+1}}\}$, respectively.
Thus, the claim follows from~\eqref{040514a}, \eqref{040514b} and~\eqref{040514e}.\medskip

Let us continue the proof of Lemma~\ref{040914b}.
Now, since, obviously, $\cA\ge\cA_i$, we conclude by \eqref{040514c} that $\cA\ge\cyc(K_i,D)$. To verify
the converse inclusion, we have to prove that every $K_i$-orbit  $Z'$ belongs to~$\cS(\cA)$. Suppose
on the contrary that some $Z'$ properly contains a set $Z\in\cS(\cA)$. Then 
\qtnl{010314a}
Z\subseteq D_{n-1}\setminus D_1.
\eqtn
Indeed, from~\eqref{010514a} and statement~(3) of Lemma~\ref{010514b}, it follows that the $K_i$-orbits
outside $D_{n-1}$ are the basic sets of~$\cA$. Besides, the orbits $\{e\}$ and $\{c_1\}$ are also basic sets,
because $\cA\ge\cyc(K_i,D)$. Finally, the orbit $\{s,c_1s\}$ belongs to $\cS(\cA)$ by the second 
part of statement~(1) of Lemma~\ref{010514b}.\medskip

From \eqref{010314a} and the first part of statement~(1) of Lemma~\ref{010514b}, it follows that the set $Z$ is regular.
So it is an orbit of an automorphism group of $C$. This implies that $Z$ has cardinality $1$, $2$, or $4$.
The latter case  is impossible, because otherwise $Z=Z'$. We claim that the first case is also impossible. Indeed, 
otherwise $Z=\{z\}$ for some $z\in D_{n-1}\setminus D_1$. Then 
$zX$ is a highest basic set of~$\cA$. However, $(zX)_0=\{zx,z\varepsilon x^{-1}\}$.
Therefore, $\varepsilon (zx)^{-1}=z\varepsilon x^{-1}$, and hence $z=z^{-1}$. Thus, $z\in D_1$. Contradiction.\medskip

To complete the proof  of Lemma~\ref{040914b}, let $|Z|=2$. Without loss of generality, we
can assume that the order of an element in $Z$ 
is minimal possible. Clearly, $|Z'|=4$, and hence 
$$
Z'=\{z^{\pm 1},c_1z^{\pm 1}\},
$$
where either $z=x^2$, or $z=sx^{2^m}$ with $m\in\{3,\ldots,n-3\}$ (see Table~\ref{tblsg1}). Choose $z\in Z$ so that 
$$
Z=\{z,c_1z\}\qoq Z=\{z,\varepsilon' z^{-1}\},
$$
where $\varepsilon'\in\{e,c_1\}$. In the first case, the singleton $Z^2=\{z^2\}$ is a basic set of $\cA$. By the above 
argument, this implies that $z^2\in D_1$. Thus, $z\in D_2$. Contradiction. 
In the second case, $\cA$ contains the element 
$$
(sz+sz^{-1})(z+\varepsilon' z^{-1})=sz^2+s\varepsilon' z^{-2}+s+s\varepsilon'.
$$
Since $s\not\in\cA$, this implies that $\varepsilon'=c_1$. Therefore, $\{sz^2,s\varepsilon' z^{-2}\}$ cannot be a basic
set of $\cA$, and hence $m\ne 3$. But then, the minimality of $|z|$ implies that 
$\{sz^2,s\varepsilon' z^{-2},s,s\varepsilon'\}\in\cS(\cA)$ that is impossible, because $s+sc_1\in\cA$.\bull

\lmml{040914c}
Let $|X_0|=1$  and $X_1\ne\emptyset$. Then $\cA=\cyc(K_i,D)$
with $i\in\{8,9,10,11\}$.
\elmm
\proof In this case $X=\{x,ys\}$, where $\grp{y}=C$. It follows that $X$ generates $D$. Moreover, $n\ge 3$, because 
$c_1\not\in\rad(X)$. Therefore, there are two orthogonals in the S-ring~$\cA$: one of them is in  
$(X\,X^{-1})\cap Cs$, and another one is in $(X\,X^{-1})\cap C$. Thus, $\cA$ is quasi-thin
 by Lemma~\ref{100914a}, and hence  schurian by Lemma~\ref{100914y}. The latter implies also that
the stabilizer $K$ of the point $e$ in the group $\aut(\cA)$, has a faithful regular orbit. Therefore, the index of $D$ in $\aut(\cA)$ 
is equal to~$2$. But then $D\trianglelefteq \aut(\cA)$, and hence $K\le\aut(D)$.
Consequently, $\cA=\cyc(K,D)$ and the group $K$ is generated by an involution
$\sigma\in\aut(D)$. This involution interchanges $x$ and $ys$. Therefore, the automorphism $\sigma$ 
is uniquely determined. We leave the reader to verify that $\sigma$ is one of automorphisms that are listed 
in the rows 8, 9, 10, 11 of Table~\ref {280414b}.\bull

\section{S-rings over $D=\mZ_2\times\mZ_{2^n}$: automorphism groups in regular case}\label{231014a2}

In this section we find the automorphism group of a regular S-ring over $D$ with trivial radical. For this purpose, we 
need the following concept introduced in~\cite{KK}: a permutation group is called {\it $2$-isolated} if no other
group is $2$-equivalent to it. The following statement is the main result of this section; 
it shows, in particular, that  a regular S-ring over $D$ with trivial radical is normal.

\thrml{090514a}
Let $\cA$ be a regular S-ring over the group $D$. Suppose that $\rad(\cA)=e$. Then for any
$\cA$-group $L$ of order at most $2$, the group $\aut(\cA_{D/L})$ is $2$-isolated. In
particular, if $\cA=\cyc(K,D)$ for some $K\le\aut(D)$, then
$\aut(\cA)=DK$.
\ethrm

The proof will be given in the end of the section. The following statement provides a sufficient condition
for a permutation group to be $2$-isolated.

\lmml{160514a}
Let $\cA$ be an S-ring and $G=\aut(\cA)$. Suppose that the point stabilizer of $G$ has a faithful regular
orbit. Then the group $G$ is $2$-isolated.
\elmm
\proof  It was proved in \cite[Theorem~$3.5'$]{KK} that $G$ is 
$2$-isolated whenever it is $2$-closed and a two-point stabilizer 
of $G$ is trivial. However, the latter exactly means that a 
point stabilizer of~$G$ has a faithful regular orbit.\bull\medskip

To apply Lemma~\ref{160514a}, we need the following auxiliary statement giving a sufficient condition 
providing the existence of a faithful regular orbit of a point stabilizer in the automorphism group
of an S-ring.

\lmml{150514d}
Let $\cA$ be an S-ring over an abelian group $H$. Suppose that a set $X\in\cS(\cA)$ satisfies the
following conditions:
\nmrt
\tm{1} $\grp{\tr(X)}=H$,
\tm{2} $c_{XY}^Z=1$ for each $Z\in\cS(\cA)_{\tr(X)}$ and some $Y\in\cS(\cA)$,
\tm{3} $c_{XY}^X=1$ for each $Y\in\cS(\cA)_{XX^{-1}}$.
\enmrt
Then $X$ contains a faithful regular orbit of the group $\aut(\cA)_e$.
\elmm
\proof  Denote by $\cX$ the Cayley scheme over $H$ associated with the S-ring~$\cA$. Then the
relation $r=R_H(\tr(X))$ is a union of basic relations of $\cX$. Clearly, $r$ is symmetric, and 
connected (condition~(1)). Moreover, conditions~(2)
and~(3) imply that the coherent configuration $(\cX_e)_{er}$ is semiregular. Thus, by \cite[Theorem~3.3]{P13}
given $x\in X$, the two-element set $\{e,x\}$ is a base of the scheme $\cX$, and hence of the group $\aut(\cX)$. 
This implies that $\{x\}$ is a base of the group  $K=\aut(\cX)_e$. Thus, $x^K\subseteq X$ is a faithful regular
orbit of~$K$.\bull

{\bf Proof of Theorem~\ref{090514a}.}\ From Theorem~\ref{230214a}, it follows that $\cA=\cyc(K,D)$, where $K=K_i$ is one of the groups 
listed in Table~\ref{280414b}, $1\le i \le 11$. Therefore, taking into account that the groups $DK$ and $\aut(\cA)$ are $2$-equivalent,
we conclude that the second part of the theorem statement immediately follows from the first one. To prove the latter, without loss of generality, we assume that $i\ge 2$ and $n\ge 4$. \medskip
 
Let $L\le D$ be an $\cA$-group. In what follows, we set $H=D/L$ and $\pi=\pi_L$. To prove that the group $\aut(\cA_H)$ is $2$-isolated, 
it suffices to verify that its point stabilizer has a faithful regular orbit (Lemma~\ref{160514a}). The remaining part of the proof 
is divided into three cases.\medskip

{\bf Case 1:} $L=\grp{\varepsilon s}$,  where $\varepsilon\in\{e,c_1\}$. Here, $H$ is a cyclic group and the S-ring $\cA_H=\cyc(\pi(K),H)$
is cyclotomic. Moreover, $i\in\{2,3,8,9\}$ by statement~(3) of Corollary~\ref{100514d}. Therefore, the order of the group $\pi(K)\le\aut(H)$ 
is at most~$2$. By the implication $(3)\Rightarrow(2)$ of~\cite[Theorem~6.1]{EP02}, this implies
that the group $\aut(\cA_H)_e$ has a faithful regular orbit. Thus, the group $\aut(\cA_H)$ is $2$-isolated by Lemma~\ref{160514a}.\medskip

{\bf Case 2:} $e\le L \le\grp{c_1}$ and $|K|=2$. Here, $i\not\in\{1,6,7\}$ and each basic set  of $\cA$ is of cardinality at most~$2$.
Since $\cA$ is commutative, the latter is also true for the basic sets of~$\cA_H$. Therefore, this S-ring is quasi-thin. So by the 
second part of 
Lemma~\ref{100914y} it suffices to prove that $\cA_H$ has at least two orthogonals.  
To do this let $X$ be a basic set of $\cA$ that contains a generator of~$C$. Since
the S-ring $\cA$ is cyclotomic, $\pi(X^{(2)})$ and $\pi(X^{(4)})$ are basic sets of $\cA_H$.
Moreover, they are distinct, because $n\ge 4$. Finally, they are orthogonals,
because $\pi(X^{(2)})\subseteq \pi(X)\,\pi(X^{-1})$ and 
$\pi(X)^{(4)}\subseteq \pi(X^{(2)})\,\pi(X^{(-2)})$.\medskip

{\bf Case 3:} $e\le L \le\grp{c_1}$ and $|K|=4$.  Here $i=6$ or $7$. It suffices to verify that the hypothesis of 
Lemma~\ref{150514d} is satisfied for  a highest basic set $X$ of the S-ring~$\cA_H$. To do this we first observe 
that the sets $X_0$ and $X_1$ are not empty. Therefore,
$\tr(X)=D\setminus D_{n-1}$, and condition~(1) is, obviously, satisfied.\medskip

To verify conditions (2) and (3), suppose first that $L=e$. Then  for $x\in X_0$, we have
$$
X=\{x,\varepsilon x^{-1},sc_2x,sc_2^{-1}\varepsilon x^{-1}\},
$$
where $\varepsilon\in\{e,c_1\}$. Since $n\ge 4$, the elements $xy^{-1}$ with $y\in X$
belong to distinct $K$-orbits  of  cardinalities $1,2,2$ and $4$.  A straightforward check
shows that if $Y$ is one of these orbits, then
$$
|Y|c_{XX^{-1}}^Y=4.
$$ 
Therefore, $4=|Y|c_{XX^{-1}}^Y=|X|c_{X^{}Y^{-1}}^X$, and condition~(3) is satisfied,
because $|X|=4$.  Let now $Z\in\cS(\cA)_{\tr(X)}$. Then
$$
Z=\{xy,\varepsilon (xy)^{-1},sc_2xy,sc_2^{-1}\varepsilon (xy)^{-1}\}
$$
for some $y\in C_{n-1}$. Let $Y$ be the set $\{y^{\pm 1},c_1y^{\pm 1}\}$ , 
$\{y^{\pm 1}\}$, or $\{y\}$ depending on whether $y$ belongs to  $C_{n-1}\setminus C_{n-2}$,
$C_{n-2}\setminus C_1$, or $C_1$, respectively. Then $Y$ is a basic set of $\cA$. Moreover,
a straightforward computation shows that in any case, $c_{YZ}^X=1$. 
Since $c_{X^{}Y^{}}^{Z^{}}=c_{Y^{-1}Z^{-1}}^{X^{-1}}=c_{Y^{}Z^{}}^{X^{}}$,  
condition~(2) is also satisfied.\medskip

Let now $L=\grp{c_1}$. To simplify notations, we identify the group $H=D/L$ with
$D_{n-1}$, write $\cA$ instead of $\cA_H$,\footnote{In our case, the S-ring $\cA_H$ does
not depend on the choice of $i\in\{5,6\}$.} and use the notation $x$ and $s$ for the 
$\pi$-images of $x$ and $sc_1$, respectively. Thus, $\cA$ is a cyclotomic
S-ring over $D_{n-1}$ and 
$$
X=\{x,x^{-1},sx,sx^{-1}\}
$$
is a highest basic set of $\cA$. It follows that $C_{n-2}$ is an $\cA$-group and any basic set inside $C_{n-2}$ is of 
the form $\{z^{\pm 1}\}$ for a suitable $z\in C_{n-2}$. Since $sC_{n-2}$ is an $\cA$-set, the elements $xy^{-1}$ 
with $y\in X$ belong to distinct basic sets of~$\cA$. Therefore, the set $X^{}\,X^{-1}$ 
consists of basic sets $Y$ for which $c_{XY}^X=1$. Thus, condition~(3) is satisfied.  
Let now $Z\in\cS(\cA)_{\tr(X)}$. Then
$$
Z=\{xy,(xy)^{-1},sxy,s (xy)^{-1}\}
$$
for some $y\in C_{n-2}$. Taking $Y$ to be the basic set $\{y^{\pm 1}\}$, we find that 
$c_{YZ}^X=1$. Thus, condition~(2) is also satisfied, and we are done.\bull

\section{S-rings over $D=\mZ_2\times\mZ_{2^n}$:  non-trivial radical case}\label{140914b}

In Theorems~\ref{170414a} and~\ref{230214a}, we completely described the structure of an S-ring over~$D$ that has
trivial radical. In this section, we study the remaining S-rings.

\thrml{060414a}
Let $\cA$ be an S-ring over a group $D$. Suppose that $\rad(\cA)\ne e$. Then $\cA$ is a proper generalized
$S$-wreath product, where the section $S=U/L$ is such that
\qtnl{171114a}
\cA_S=\mZ S\qoq |S|=4\qoq \rad(\cA_U)=e\ \text{and}\ |L|=2.
\eqtn
\ethrm
\proof Denote by $U$ the subgroup of~$D$ that is generated by all $X\in\cS(\cA)$ such that $\rad(X)=e$.

\lmml{170414b}
$U$ is an $\cA$-group and $\rad(\cA_U)=e$.
\elmm
\proof The first statement is clear. To prove the second one, without loss of generality, we 
can assume that $U=D$. Then there exists a highest set $X\in\cS(\cA)$ such that
$\rad(X)=e$. Suppose first that $X$ is not regular. Then $X\cap E\ne\emptyset$ by 
Theorem~\ref{290714a}. Therefore, $X$ is one of basic sets $X_{c_1}$, $X_s$, or $X_{sc_1}$
from Theorem~\ref{250814a}. Since $\rad(X)=e$, we have $X=X_{c_1}$ and $X=X_s$
in statements~(2) and~(3) of this theorem,  respectively. Moreover, since $X$ is highest, $D=\grp{X}$
in statements~(1), (2), (3). Therefore, in these three cases $\rk(\cA)=2$, and hence $\rad(\cA)=e$. In the 
remaining case (statement~(4) of Theorem~\ref{250814a}), we have $D=H$, and hence
$\cA=\cA_U\otimes\cA_{\grp{s}}$. Since each of the factors is of rank~$2$, this
implies that again  $\rad(\cA)=e$.\medskip

From now on,  we can assume that any highest basic set of~$\cA$ with trivial radical, 
is regular (Corollary~\ref{011213d}). Moreover, if $X$ is one of them  and both $X_0$ and
$X_1$ are not empty, then all highest basic sets are pairwise rationally conjugate and, hence
$\rad(\cA)=\rad(X)=e$. Thus, we can also assume that $\grp{X}$ is a cyclic group $C$
of order at least~$4$ that has index $2$ in $D$.\medskip

Since $\rad(\cA_C)=\rad(X)=e$, the circulant S-ring $\cA_C$ is cyclotomic (see 
Subsection~\ref{141114b}). Together with $|C|\ge 4$, this shows that $c_1\in\cA$. 
We claim that any $Y\in\cS(\cA)$ such that $\rad(Y)=e$ and $D=\grp{X,Y}$,
is regular.  Indeed, otherwise by Theorem~\ref{290714a}, we have $Y=X_h$,
where $h$ is a non-identity element of the group~$E$. However, $h\ne c_1$: otherwise $Y=\{c_1\}$ by above,
and $\grp{X,Y}=C$, in contrast to the assumption. Since $X_{c_1}=\{c_1\}$ and $\rad(Y)=e$, 
only statement~(4) of Theorem~\ref{250814a} can hold.  But then, $Y$ is a singleton 
in $E$, and hence it is regular. Contradiction.\medskip

To complete the proof, let $Y$ be a regular basic set with trivial radical. We can assume
that $Y\subseteq D\setminus C$, for otherwise $D=C$ and $\rad(\cA)=\rad(X)=e$.
If $Y$ is not highest, then any basic set $Z\subseteq XY$ is highest, and $\grp{X,Z}=D$. 
Moreover, $\rad(Z)=e$ by Lemma~\ref{150414a}. Thus, we can also assume that $Y$ is highest. 
Then any highest basic set of $\cA$ is rationally conjugate to either $X$ or $Y$. Thus,
$\rad(\cA)=e$, as required.\bull

By the theorem hypothesis and Lemma~\ref{170414b}, we have $U\ne D$. We observe also that
by Theorem~\ref{030414a},  the group $U$ contains every minimal $\cA$-group.

\lmml{130914a}
Suppose that there is a unique minimal $\cA$-group, or $c_1\in\rad(X)$ for all $X\in\cS(\cA)_{D\setminus U}$. Then the statement of Theorem~\ref{060414a} holds.
\elmm
\proof Let $L$ be a unique minimal $\cA$-group. Then the definition of $U$ implies that 
$\cA$ is a proper generalized $S$-wreath
product, where $S=U/L$. If, in addition, $|L|\le 2$, then the third statement in~\eqref{171114a} follows
from Lemma~\ref{170414b} and we are done. Let now $|L|>2$. Then $\grp{c_1}$ is not
an $\cA$-group. By statement~(1) of Corollary~\ref{100514d}, this implies that the S-ring 
$\cA_U$ is  not regular. So, by  the first part of Theorem~\ref{170414a}, it is rational.
Now, by the second part of this theorem, the uniqueness of~$L$ implies that $U=L$. Thus, 
$|S|=1$ and the first statement in~\eqref{171114a} holds.\medskip

To complete the proof, suppose that there are at least two minimal $\cA$-groups. Then, obviously,
one of them, say $H$, contains $c_1$. Therefore, $c_1\in H\le U$. On the other hand, by the 
lemma hypothesis, $c_1\in\rad(X)$ for all $X\in\cS(\cA)_{D\setminus U}$. Thus, $\cA$ is a proper generalized
$S$-wreath product, where $S=U/H$. Without loss of generality, we can assume that $|H|>2$.
If the S-ring $\cA_U$ is rational, then from the second part of
Theorem~\ref{170414a}, it follows that there is another minimal $\cA$-group $L$ of order $2$
and such that
$$
\cA_U=\cA_H\otimes\cA_L.
$$ 
Thus, $|S|=|L|=2$, and the first statement in~\eqref{171114a} holds. In the remaining
case, $\cA_U$ is a regular S-ring  by the first part of Theorem~\ref{170414a}. By statement~(1)
of Corollary~\ref{100514d}, this implies that $H=\grp{c_1}$.  Thus, $|H|=2$, and
the third statement in~\eqref{171114a} holds.\bull

Denote by $V$ the union of all sets $X\in\cS(\cA)$ such that $\rad(X)=e$ or $c_1\in\rad(X)$.
Then, obviously, $U\subseteq V$ and $V$ is an $\cA$-set. By Lemma~\ref{130914a}, we can
assume that $V\ne D$, and that $U$ contains two distinct minimal $\cA$-groups. It is easily seen that in this case
$E\subseteq V$.

\lmml{130914c}
In the above assumptions let $X\in\cS(\cA)_{D\setminus V}$. Then
\nmrt
\tm{1} $\rad(X)=\grp{s}$ or $\grp{sc_1}$,
\tm{2} $X$  is a regular set such that both $X_0$ and $X_1$ are not empty.
\enmrt
\elmm
\proof Since $X\not\subseteq U$, we have $\rad(X)\ne e$. Besides, $c_1\not\in\rad(X)$ by the
definition of $V$. Thus, statement (1) holds, because $\grp{s}$ and $\grp{sc_1}$ are the
only subgroups of $D$ that do not contain $c_1$. To prove statement (2), set
$$
L=\rad(X)\qaq\pi=\pi_{D/L}.
$$ 
Then $\rad(\pi(X))=e$. However, $D/L$ is a cyclic $2$-group by statement~(1). Therefore, 
$\pi(X)$ is the basic set  of a circulant S-ring $\cA_{D/L}$. From the description
of basic sets of such an S-ring given in Subsection~\ref{141114b}, it follows
that $\pi(X)$ is regular or is of the form
$$
\pi(X)=\pi(H)^\#
$$ 
for some $\cA$-group $H\ge D_1$ such that $|H/L|\ge 4$.  In the latter case,
$X=H\setminus L$, and hence $L$ is a unique minimal $\cA$-group in contrast to 
our assumption on $U$. Thus, the set $\pi(X)$ is regular. This implies that the
set $X$ is also regular. Finally, the fact that $X_0$ and $X_1$ are not empty,
immediately follows from statement~(1).\bull

By statement~(2) of Lemma~\ref{130914c},  the union of $\tr(X)$, where $X$ runs over the set 
$\cS(\cA)_{D\setminus V}$, is of the form $D\setminus D_k$ for some $k\ge 1$. However, 
the set $V$ coincides with the complement to this union. Thus, $V=D_k$ is an $\cA$-group.
A similar argument shows that $D_m$ is an $\cA$-group for all $m\ge k$. 

\lmml{130914u}
Let $m=\max\{2,k\}$. Then the group $L:=\rad(X)$ does not depend on the choice of
$X\in\cS(\cA)_{D\setminus D_m}$.
\elmm
\proof Suppose on the contrary that there exist basic sets $X$ and $Y$ outside $D_m$ such that
$\grp{Y}\subsetneq\grp{X}$ and $\rad(X)\ne\rad(Y)$. Then by statement (1) of Lemma~\ref{130914c}, 
without loss of generality, we can assume that 
\qtnl{181114a}
\rad(X)=\grp{s}\qaq \rad(Y)=\grp{sc_1}.
\eqtn
By statement~(2) of that lemma, $X_0$ is a regular non-empty set. Therefore, it
is an orbit of a subgroup of $\aut(C)$. Moreover, from the first equality in \eqref{181114a}, 
it follows that  $\rad(X_0)=e$. Let now  $\pi:D\to D/\grp{s}$ be the quotient epimorphism. Then 
$\pi(X)=X_0$ is a basic set of a circulant S-ring $\cA'=\pi(\cA)$. It follows that
$\cA'_{\grp{X_0}}$ is a cyclotomic S-ring with trivial radical. Moreover, $Y'=\pi(Y)$ is a 
basic set of this S-ring and $|\grp{Y'}|\ge 2^{m+1}\ge 8$. Therefore, by Lemma~\ref{021109a} applied 
for $\cA=\cA'_{\grp{X_0}}$ and $S=\grp{Y'}/e$, we obtain that
$$
\rad(\cA'_{\grp{Y'}})=e.
$$
This implies that $\rad(Y')=e$. On the other hand, by the second equality in~\eqref{181114a},
we have $\rad(Y')=\grp{c_1}\ne e$. Contradiction.\bull

By Lemma~\ref{130914u}, the S-ring $\cA$ is the generalized $S$-wreath product,
where $S=D_2/L$ if $k=1$, and $S=V/L$ if $k\ge 2$. The only case when this generalized
product is not proper, is $D=D_2$ and $k=1$. However, in this case $\cA$ is, obviously, a proper
$E/L$-wreath product and $|E/L|=2$. Thus, if $k\le 2$, then the first or the second statement 
in~\eqref{171114a} holds, and we are done. To complete the proof of Theorem~\ref{060414a},
it suffices to verify that the third statement in~\eqref{171114a} holds whenever $k\ge 3$.
But this immediately follows from the lemma below. 

\lmml{220414u}
If $k\ge 3$, then $\rad(\cA_V)=e$. In particular, $V=U$.
\elmm
\proof The second statement follows from the first one and statement~(1) of Lemma~\ref{130914c}.
To prove that $\rad(\cA_V)=e$, let $X$ be a highest basic set of the S-ring $\cA_V$. Since $V\ne D$, there exists 
a set $Y\in\cS(\cA)_{D\setminus V}$ such that $X\subseteq Y^2$. By Lemma~\ref{130914c}, 
we have  $Y=LY_0$, where $Y_0$ is an orbit of a subgroup of $\aut(C)$ such that  $\rad(Y_0)=e$. 
However, $Y_0=\{y\}$ or $Y_0=\{y,\varepsilon y^{-1}\}$, where $\varepsilon\in\{e,c_1\}$
(statement~(1) of Lemma~\ref{211113a}). Therefore,
\qtnl{240414b}
Y^2=LY_0^2=L\times\css
\{y^2\},                       &\text{if $Y_0=\{y\}$},\\
\{\varepsilon,y^{\pm 2}\},     &\text{if $Y_0=\{y,\varepsilon y^{-1}\}$}.\\
\ecss
\eqtn
On the other hand, since $X\subset V$, the definition of $V$ implies that $\rad(X)=e$ or $c_1\in\rad(X)$. In 
the former case, $\rad(\cA_V)=e$, and we are done. Suppose that $c_1X=X$. Then the set $Y^2$ contains
a $\grp{c_1}$-coset $\{x,c_1x\}$ for all $x\in X$. By~\eqref{240414b}, this implies that
$$
\{x,c_1x\}\subset \{\varepsilon,y^{\pm 2}\}
$$
which is impossible, because $|x|=2^k\ge 8$.\bull

\section{S-rings over $D=\mZ_2\times\mZ_{2^n}$: schurity}\label{231014a3}

In this section, based on the results obtained in Sections~\ref{140914a}--\ref{140914b}, we prove the following main 
theorem.

\thrml{090514b}
For any integer $n\ge 1$, every S-ring over the group $D=\mZ_2\times\mZ_{2^n}$ is schurian. In particular, 
$D$ is a Schur group.
\ethrm
\proof 
The induction on $n$. An exhaustive computer search of all S-rings over small groups shows that 
$D$ is a Schur group for $n\le 4$. Let $n\ge 5$. We have to verify that any S-ring $\cA$ 
over~$D$ is schurian. However, if $\rad(\cA)=e$, then this is true by Theorems~\ref{170414a}
and~\ref{230214a}. For the rest of the proof, we need the following result from~\cite{EP12}
giving a sufficient condition for a generalized wreath product of S-rings to be schurian.

\thrml{140914f}{\rm \cite[Corollary~5.7]{EP12}}
Let $\cA$ be an S-ring over an abelian group~$D$. Suppose that $\cA$ is the generalized $S$-wreath 
product of schurian S-rings $\cA_{D/L}$ and $\cA_U$, where~$S=U/L$. Then $\cA$  is schurian if and only if there exist two groups 
$\Delta_0\ge(D/L)_{right}$ and $\Delta_1\ge U_{right}$, such that
\qtnl{140914r}
\Delta_0\twoe\aut(\cA_{D/L})\qaq
\Delta_1\twoe\aut(\cA_U)\qaq 
(\Delta_0)^{U/L}=(\Delta_1)^{U/L}.
\eqtn
\ethrm

\crllrl{110514a}
Under the hypothesis of Theorem~\ref{140914f}, the S-ring $\cA$  is schurian whenever the 
group $\aut(\cA_S)$ is $2$-isolated.
\ecrllr
\proof Set $\Delta_0=\aut(\cA_{D/L})$ and $\Delta_1=\aut(\cA_U)$. Then the first two
equalities in~\eqref{140914r} hold, because the S-rings $\cA_{D/L}$ 
and $\cA_U$ are schurian. Since the group $\aut(\cA_S)$ is $2$-isolated, we have
$
(\Delta_0)^S=\aut(\cA_S)=(\Delta_1)^S,
$ 
which proves the third equality in~\eqref{140914r}. Thus, $\cA$ is schurian
by Theorem~\ref{140914f}.\bull

Let us turn to the proof of Theorem~\ref{090514b}. Now, we can assume that 
$\rad(\cA)\ne e$. Then by Theorem~\ref{060414a},  the S-ring~$\cA$  is a proper generalized $S$-wreath
product, where the section $S=U/L$ is such that formula~\eqref{171114a} holds. Besides, by induction,
the S-rings $\cA_{D/L}$ and $\cA_U$ are schurian. Suppose that $\cA_S=\mZ S$, or $|S|=4$,
or $|L|=2$ and $\cA_U$ is a regular S-ring with trivial radical. Then the group $\aut(\cA_S)$ 
is $2$-isolated: this is obvious in the first two cases and follows from Theorem~\ref{090514a}
(applied for $\cA=\cA_U$) in the third one.  Thus, $\cA$ is schurian
by  Corollary~\ref{110514a}.\medskip

To complete the proof, we can assume that $|S|=|2^m|$, where $m\ge 3$, and that
$\cA_U$ is a non-regular S-ring with trivial radical. Then 
$\cA_U=\cA_H\otimes\cA_L$, where $|H|\ge 4$ and $\rk(\cA_H)=2$ (Theorem~\ref{170414a}). 
Therefore,
\qtnl{140914w}
\aut(\cA_U)^S=(\sym(H)\times\sym(L))^{U/L}=\sym(S).
\eqtn
On the other hand, $L=\grp{s}$ or $L=\grp{sc_1}$, because $c_1\in H$. Therefore,
the S-ring $\cA_{D/L}$ is circulant. Besides, $S$ is an $\cA_{D/L}$-section of 
composite order. By \cite[Theorem~4.6]{EP12}, this implies that
\qtnl{150914a}
\aut(\cA_{D/L})^S=\sym(S).
\eqtn
By \eqref{140914w} and \eqref{150914a}, relations \eqref{140914r} are true
for the groups $\Delta_0:=\aut(\cA_{D/L})$ and $\Delta_1:=\aut(\cA_U)$. Thus, the 
S-ring~$\cA$ is schurian by Theorem~\ref{140914f}.\bull

\section{A non-schurian S-ring over $M_{2^n}$}\label{231014a4}

The main result of this section is the following theorem in the proof of which
we construct a non-schurian S-ring over the group $M_{2^n}$
defined in~\eqref{170914a}.

\thrml{071113a}
For any $n\ge 4$, the group $M_{2^n}$ is not Schur.
\ethrm
\proof The group $M_{16}$ is not Schur \cite[Lemma~3.1]{PV}. Suppose that $n\ge 5$.
Denote by $e$ the identity of the group~$G=M_{2^n}$, and
by $H$ the
normal subgroup of $G$ that is generated by the elements $c=a^{2^{n-3}}$ and $b$. Then
$H\simeq \mZ_4\times \mZ_2$ and
\qtnl{110413b}
H=Z_0\,\cup\, Z_1\,\cup\, Z_2,
\eqtn
where the sets $Z_0=\{e\}$, $Z_1=\{c^2\}$, and $Z_2=H\setminus\grp{c^2}$ are mutually disjoint. Next, let us 
fix two other decompositions of $H$ into a disjoint union of subsets: 
$$
H=\underbrace{B\,\cup\,Bc^3}_{X_1}\ \cup\ \underbrace{Bc^2\,\cup\,Bc}_{Y_1}
\quad=\quad
\underbrace{B'\,\cup\,B'c}_{X_2}\quad\cup\quad\underbrace{B'c^2\,\cup\,B'c^3}_{Y_2},
$$
where $B$ and $B'$ are the groups of order~$2$ generated by the involutions~$b$ and $b':=c^2b$.
Then a straightforward computation shows that 
\qtnl{110413a}
Ha\cup Ha^{-1}=\underbrace{X_1a\cup X_2a^{-1}}_{Z_3}\ \cup\ \underbrace{Y_1a\cup Y_2a^{-1}}_{Z_4},
\eqtn
in particular, the sets $Z_3$ and $Z_4$ are disjoint. Moreover, $Z_3c^2=Z_4$, because
$Y_1=X_1c^2$ and $Y_2=X_2c^2$. Finally, there are exactly $m'=2^{n-3}-3$ cosets of $H$ in $G$,
other than $H$, $Ha$, and $Ha^{-1}$. Let us combine them in pairs as follows
\qtnl{110413ae}
Z_{i+3}:=Ha^i\cup Ha^{-i},\quad i=2,3,\ldots m,
\eqtn
where $m=(m'-1)/2+1$ and $Z_{m+1}=Ha^{2^{n-2}}$. 
Then the sets $Z_0,Z_1,\ldots,Z_{r-1}$ with $r=m+5$ form
a partition of the group~$G$; denote it by~$\cS$. The submodule of $\mZ G$ spanned
by the elements $\und{Z_i}$, $i=0,\ldots,r-1$, is denoted by $\cA$.

\lmml{170914s} 
The module $\cA$ is an S-ring over $G$. Moreover, $\cS(\cA)=\cS$.
\elmm
\proof 
From the above definitions, it follows that $Z_i^{-1}=Z_i^{}$ for all $i$. Thus, it suffices
to verify that given $i$ and $j$, the product $\und{Z_i}\,\und{Z_j}$ is a linear combination 
of $\und{Z_k}$, $k=0,\ldots,r-1$. However, it is easily seen that
$$
\und{Ha^i}\,\und{Ha^j}=\und{Ha^j}\,\und{Ha^i}=\und{Ha^{i+j}}=\und{a^{i+j}H}
$$
for all $i,j,k$. Therefore, the required statement holds whenever $i,j\not\in\{3,4\}$.
To complete the proof, assume that $i=3$ (the case $i=4$ is considered analogously). Then
a straightforward check shows that
\nmrt
\item[$\bullet$] $\und{Z_3}\,\und{Z_1}=\und{Z_4}$,\quad $\und{Z_3}\,\und{Z_2}=\und{Z_4}$,\quad $\und{Z_3}\,\und{Z_{r-1}}=4\und{Z_{r-2}}$,
\item[$\bullet$] $\und{Z_3}\,\und{Z_3}=8\,\und{Z_0}+2\,\und{Z_5}+4\,\und{Z_2}$,
\item[$\bullet$] $\und{Z_3}\,\und{Z_4}=8\,\und{Z_1}+2\,\und{Z_5}+4\,\und{Z_2}$,
\item[$\bullet$] $\und{Z_3}\,\und{Z_5}=4\und{Z_3}+4\und{Z_4}+4\und{Z_6}$,
\item[$\bullet$] $\und{Z_3}\,\und{Z_i}=4\und{Z_{i-1}}+4\und{Z_{i+1}}$, $i=6,\ldots r-2$.
\enmrt
Since $\und{Z_3}$ commutes with $\und{Z_j}$ for all $j$, we are done. 
\bull 

By Lemma~\ref{170914s}, the statement of Theorem~\ref{071113a} immediately
follows from the lemma below.

\lmml{231112a}
The S-ring $\cA$ is not schurian.
\elmm
\proof Suppose on the contrary that $\cA$ is schurian. Then it
is the S-ring associated with the group
$\Gamma=\aut(\cA)$.  It follows that the basic set $Z_2$ is an orbit of the one-point 
stabilizer of~$\Gamma$. Since $|Z_2|=6$, there exists an element
$\gamma\in\Gamma$ such that
\qtnl{231112b}
|\gamma^{Z_2}|=3.
\eqtn
On the other hand, due to~\eqref{110413a}, the quotient S-ring $\cA_{G/H}$ is isomorphic
to the S-ring associated with the dihedral group of order $2^{n-2}$ in its natural permutation representation of degree $2^{n-3}$. 
Therefore, $\aut(\cA_{G/H})$ is a $2$-group. It contains a subgroup $\Gamma^{G/H}$,
and hence the element $\gamma^{G/H}$. So by~\eqref{231112b}, the 
permutation~$\gamma$ leaves each $H$-coset fixed (as a set). Therefore,
\qtnl{281112a}
\gamma^{H\cup Ha}\in\aut(C),
\eqtn
where $C$ is the bipartite graph with vertex set $H\cup Ha$ and the edges $(h,hax)$ with $x\in X_1$. However, 
the graph $C$ is isomorphic to the lexicographic product of the empty graph with~$2$ vertices and the
undirected cycle of length~$8$. Therefore, $\aut(C)$ is a $2$-group. By~\eqref{281112a}, this implies
that $|\gamma^H|$ is a power of~$2$. But this contradicts to~\eqref{231112b}, because $Z_2\subset H$.\bull

\section{S-rings over $D=D_{2n}$: divisible difference sets}\label{231014a6}

\sbsnt{Preliminaries.} In the rest of the paper, we deal with S-rings over a
dihedral $2$-group $D=D_{2n}$ of order $2n$. Interesting examples of such rings arise from difference sets. To
construct them, let us recall some definitions from~\cite{P95}.\medskip

 Let $T$ be a $k$-subset of a  
group $G$  of order  $mn$  such  that  every element  outside  a  subgroup  $N$  of order  
$n$  has  exactly  $\lambda_2$  representations  as  a quotient  $gh^{-1}$  with elements 
$g,h\in G$, and elements  in  $N$  different from the identity  have  exactly  
$\lambda_1$  such  representations,
\qtnl{161014a}
\underline{T}\cdot\underline{T}^{-1}= 
k\cdot e + \lambda_1\underline{N\setminus e}+\lambda_2\underline{G\setminus N}.
\eqtn
Then $T$ is  called  an {\it $(m,n,k,\lambda_1,\lambda_2)$-divisible  difference set}  in  $G$
relative to $N$.  If $\lambda_1 =0$ (resp. $n = 1$), then we say that $T$ 
is  a  {\it relative  difference set} or {\it relative $(m,n,k,\lambda_2)$-difference set} 
(resp. difference set).  A difference set $T$
is {\it trivial} if it equals $G$, $\{x\}$ or $G\setminus\{x\}$, where $x\in G$.
 
\thrml{231014a}
Let $C$ be a cyclic $2$-group. Then
\nmrt
\tm{1} any difference set in $C$ is trivial,
\tm{2} there is no  relative $(2^a,2,2^a,2^{a-1})$-difference set in $C$.
\enmrt
\ethrm
\proof Statement (1) follows from \cite[Theorem~II.3.17]{BJL} and \cite[Theorem~1.2]{DS94}.
Statement~(2) follows from \cite[Theorems~4.1.4,4.1.5]{P95}.\bull

\sbsnt{Constructions.}\label{231014w} Let $D$ be a dihedral group of order $2n$, $C$ the cyclic subgroup 
of $D$ of order $n$ and $H$ a subgroup of $C$. Let $T$ be a non-empty subset of~$C$
such that $|T\cap xH|$ does not depend on $x\in T$ ({\it the intersection condition},
cf. Lemma~\ref{090608a}). Set 
$$
\cS:=\{e,\ H\setminus e,\ C\setminus H,\ Ts,\ T's\},
$$
where $T'=C\setminus T$. Clearly, $\cS$ is a partition of $D$ such that condition (S1) 
is satisfied. Since all the elements of $\cS$ are symmetric, condition (S2) is also
satisfied. Set 
\qtnl{161014c}
\cA:=\cA(T,C)=\Span\{\und{X}:\ X\in\cS\}.
\eqtn
\thrml{p0}
In the above notation, $\cA$ is an S-ring over $D$ with $\cS(\cA)=\cS$ if and only if $T$ 
is a divisible difference set in $C$ relative to $H$.
\ethrm
\proof To prove the ``only if'' part, suppose that $\cA$ is an S-ring with $\cS(\cA)=\cS$.  
Then $Ts$ is a basic set of $\cA$ and $Ts\,Ts=T\,T^{-1}$ is a subset of~$C$. 
Therefore,
\qtnl{191014a}
\underline{T}\cdot\underline{T}^{-1}=
\underline{Ts}^2=
|T|e+\lambda_1\underline{H\setminus e} + \lambda_2\,\underline{C\setminus H},
\eqtn
where $\lambda_1=c_{Ts\,Ts}^{H\setminus e}$ and 
$\lambda_2=c_{Ts\,Ts}^{C\setminus H}$. Thus, $T$ is a divisible difference 
set in $C$ relative to~$H$.\medskip

To prove the ``if'' part, suppose that $T$ is a divisible difference set in $C$ relative
to~$H$. It suffices to verify that $\cA\cdot\cA\subseteq \cA$. To do this, denote by $\cA'$
the module spanned by $e$, $\underline{H}$, $\underline{C}$, and $\underline{D}$. 
Then, obviously, $\cA'\cdot\cA'\subseteq\cA'$ and $\cA =\Span\{\cA',\underline{sT}\}$. 
Thus, we have to check that 
$$
\underline{sT}^2\in\cA\qaq \cA'\cdot\underline{sT}\,\subseteq\, \cA.
$$
The first inclusion follows from~\eqref{161014a}, because $T$ is a divisible difference set. 
Routine calculations show that the second inclusion is equivalent to the inclusion
$\und{T}\cdot\und{H}\in\cA$. However, this easily follows from the intersection condition.\bull

We do not know any divisible difference set over a cyclic $2$-group that satisfies the intersection condition.
However, we can slightly modify the construction by taking the set $T$ to satisfy the intersection condition, but
this time inside $C\setminus H$. Then using the same argument, one can
construct an S-ring of rank $6$ that coincides with $\cA$ on $C$ and has three
basic sets outside $C$: $Ts$, $T's$, and $Hs$, where $T'=C\setminus (T\cup H)$.\medskip 

The S-rings of this form do exist. It suffices to take a classical relative 
$(q+1,2,q,(q-1)/2)$-difference set $T$ defined as follows (see~\cite[Theorem~2.2.13]{P95}). 
Take an  affine line  $L$ in a  $2$-dimensional linear space over finite field $\mF_q$ that does not contain  
the origin. Then $L$ is a relative $(q+1,q-1,q,1)$-difference set in the 
multiplicative group of $\mF_{q^2}$. Let $\pi$ be a quotient epimorphism from this group
onto $C$ such that $|\ker(\pi)|=m$ for some divisor $m$ of $q-1$. Then $\pi(L)$ is a 
cyclic relative  $(q+1,(q-1)/m,q,m)$-difference set in~$C$. When $C$ is a $2$-group, the 
number $q+1$ is a $2$-power, and so, $q$ is a Mersenne number. If, in addition,
$m=(q-1)/2$, then we come to the required set $T$.

\crllrl{120614a}
Let $\cA$ be an S-ring over a dihedral $2$-group $D$. Suppose that $C$ is an $\cA$-group. Then
\nmrt
\tm{1} if $\rk(\cA_C)=2$, then $\cA\cong\cA_C\ast\cB$, where $\ast\in\{\wr,\cdot\}$ 
and $\cB=\mZ\mZ_2$,
\tm{2} if $\rk(\cA_C)=3$ and $\rk(\cA)=5$, then $\cA=\cA(T,C)$, where $T$ is 
a divisible difference set in $C$. 
\enmrt
\ecrllr
\proof To prove statement~(1), suppose that $\rk(\cA_C)=2$. Let $X$ be a basic set outside $C$.
Then  $X=Ts$ for some set $T\subseteq C$. From~\eqref{191014a} with $H=e$, it follows that $T$ is 
a difference set in~$C$. By statement~(1) of Theorem~\ref{231014a}, this implies that either 
$T=C$, or $T$ or $C\setminus T$ is a singleton. It 
is easily seen, that $\cA$ is isomorphic to $\cA_C\wr\cB$ in the former case, and to $\cA_C\otimes \cB$
in the other two.\medskip

To prove statement (2), suppose that $\rk(\cA_C)=3$ and $\rk(\cA)=5$. Then 
$$
\cA_C=\Span\{e,\und{H},\und{C}\}
$$ 
for some $\cA$-group $H<C$. Besides, any $X\in\cS(\cA)_{D\setminus C}$ is of the form
$X=Ts$ for some set $T\subseteq C$ satisfying the intersection condition (Lemma~~\ref{090608a}).
Thus, $\cA=\cA(T,C)$. So, $T$ is a divisible difference set in $C$ by Theorem~\ref{p0}.\bull

\sbsnt{Schurity.}
The main goal of this subsection is to prove the following theorem showing that the first
construction given in the previous subsection, produces mainly non-schurian S-rings.

\thrml{081014w}
Let $T$ be a divisible difference set in a cyclic $2$-group~$C$ relative to 
a group $H\le C$. Suppose that the intersection condition holds and
$HT\ne T$. Then the S-ring $\cA(T,C)$ defined in~\eqref{161014c} is not schurian.
\ethrm

We will deduce this theorem in the end of this subsection from a general statement on 
schurian S-rings over a dihedral $2$-group. This statement shows that if $T$ is a 
divisible difference set in a cyclic $2$-group relative to a subgroup $H$, and
the S-ring of rank $6$ associated with $T$, is schurian but not a proper generalized 
wreath product, then $H$ is of order~$2$.

\thrml{231014s}
Let $\cA$ be a schurian S-ring over a dihedral $2$-group $D$ and $H<C$ a minimal 
$\cA$-group. Suppose that $\cA$ is not a proper generalized wreath product. Then $|H|=2$.
\ethrm
\proof By the hypothesis, $\cS(\cA)=\orb(G_e,D)$ for some group $G\leq\sym(D)$ containing $D_{right}$. 
Since $H$ is an $\cA$-group, the partition 
$D/H$ of the group $D$ into the right $H$-cosets,  forms an imprimitivity system for~$G$. Denote by $N$ 
the stabilizer of this partition in $G$, 
$$
N=\{g\in G:\ (Hx)^g=Hx\ \,\text{for all}\ \,x\in D\}.
$$
Since $H_{right}\leq N$, the group $N^X$ is transitive for each block $X\in D/H$.
The following statement can also be deduced from~\cite[Lemma 2.1]{KMM11}.

\lmml{p1} 
For each block $X\in D/H$, the group $N^X$ is $2$-transitive.
\elmm
\proof Let $X\in D/H$. Then the group $G^X$ is $2$-transitive, because $\rk(\cA_H)=2$. Moreover, it contains a regular cyclic 
subgroup isomorphic to $H_{right}$. Since $H$ is a $2$-group, 
the classification of primitive groups having a regular cyclic subgroup~\cite{Jo02},
implies that $G^X$ contains 
a unique minimal normal subgroup $K$ that is also $2$-transitive.  However,  $N^X$ is a non-trivial normal subgroup 
of $G^X$. Therefore, $N^X$ contains $K$, and hence is $2$-transitive.\bull\medskip

Let us define an equivalence relation $\sim$ on the $H$-cosets by setting $X\sim Y$ if and 
only if the actions of $N$ on~$X$ and $Y$ have the same permutation character. Then by 
Lemma~\ref{p1} and a remark in~\cite[p.2]{C72}, the group $N_e$ acts transitively on each 
$X$ not equivalent to $H$. Denote by $U$ the class of $\sim$ that
contains $H$. Then each orbit of $G_e$ outside $U$ is a union of $N$-orbits. Thus,
$\cA$ is the generalized $U/H$-wreath product. By the theorem hypothesis, this implies that $U=D$.
Therefore, all classes of the equivalence relation $\sim$ are singletons. So by Lemma~\ref{p1},
we have
$$
|\orb(N,X\times Y)|=2\quad\text{for all}\ X,Y\in D/H.
$$
This yields 
us two symmetric block-designs between $X$ and $Y$ which are complementary to each other.
Since $N$ contains a cyclic subgroup $H$ which acts regularly on $X$ and $Y$, these 
block-designs are circulant, and so correspond to cyclic difference sets. 
By statement~(1) of Lemma~\ref{231014a}, they are trivial. Thus, the group 
$N_x$ with $x\in X$, has two orbits on $Y$ of cardinalities $1$ and $|Y|-1$.\medskip

To complete the proof, suppose on the contrary that $|H|>2$. Then, obviously, $|Y| > 2$. Therefore, 
the group  $N_e$ fixes exactly one point in each $H$-coset~$Y$. This implies that the set $F$ of
all fixed points of $N_e$, is of cardinality $[D:H]$. On the other hand, by 
\cite[Proposition~5.2]{KMM11},  the set $F$ is a block of $G$. Therefore, $F$ is a subgroup of~$D$. 
Moreover, since $F\cap H=e$, it is a complement for $H$ in $D$. Thus, $H=C$. Contradiction.\bull

{\bf Proof of Theorem~\ref{081014w}.} Suppose on the contrary that the S-ring $\cA=\cA(T,C)$ 
is schurian. Then the hypothesis of Theorem~\ref{231014s} is satisfied, because
$H$ is a minimal $\cA$-group and $HT\ne T$.
Thus, $|H|=2$. Denote by $x$ the element of order $2$ in~$H$. Then $x\in\cA$, and hence
$x(Ts)=T's$. This implies that $|Ts|=|T's|=m$, where $m=|C|/2$, and that $x$ appears 
neither in $\und{Ts}^2$, nor in $\und{T's}^2$. Therefore, $T$ has parameters
$(m,2,m,0,m/2)$. It follows that $T$ is a relative $(m,2,m,m/2)$-difference set in $C$.
However, this contradicts part~(2) of Theorem~\ref{231014a}.\bull

\section{S-rings over $D=D_{2^{n+1}}$: a unique minimal $\cA$-group not in $C$}

In this section we deal with S-rings over a dihedral group $D=D_{2^{n+1}}$ of
order $2^{n+1}$
and keep the notation of Subsection~\ref{011014b}. The main result here is given
by the following statement.

\thrml{011014a}
Let $\cA$ be an S-ring over the dihedral group $D$. Suppose that there
is a unique minimal $\cA$-group $H$, and that $H\not\le C$.
Then $\cA$ is isomorphic to an S-ring over $\mZ_{2^{}}\times\mZ_{2^n}$.
\ethrm
\proof The hypothesis on $H$ implies that every basic set $X$ outside 
$H$ is mixed, for otherwise $\grp{X}$ contains a non-identity 
$\cA$-subgroup of $C$. Moreover, either $H=\grp{s}$ for some $s\in D\setminus C$,
or $H$ is a dihedral group. Let us consider these two cases 
separately.\medskip

{\bf Case 1: $H=\grp{s}$} for some $s$. In this case, all basic sets except
for $\{e\}$ and $\{s\}$ are mixed.  By statement~(1) of Lemma~\ref{240913a}, this
implies that $X_0^{-1}=X_0^{}$ for all $X\in\cS(\cA)$. Besides,
$Xs\in\cS(\cA)$, because $s\in\cA$, and
$(Xs)_0=X_1$ and $(Xs)_1=X_0$. 
Thus, $X_1^{-1}=X_1^{}$ also for all $X$.\medskip

Denote by $\sigma$ the automorphism of $D$ that takes $(c,s)$ to $(c^{-1},s)$,
where $c$ is a generator of $C$. Then by the above paragraph, we have
$$
X^\sigma=(X_0\,\cup\,X_1s)^\sigma=X_0^{-1}\,\cup\,X_1^{-1}=
X_0\,\cup\,X_1s=X
$$
for all $X\in\cS(\cA)$. Therefore, the semidirect product
$D\rtimes\grp{\sigma}\le\sym(D)$ is an automorphism group of~$\cA$. 
The element $s\sigma$ of this group has order two and commutes with~$c$. Therefore, the group $D'=\grp{s\sigma,c}$ 
is isomorphic to $\mZ_2\times\mZ_{2^n}$. On the other hand, $D'$ is a regular subgroup
in $\sym(D)$. Thus, the Cayley scheme over $D$ associated with $\cA$ is
isomorphic to a Cayley scheme over $D'$. Consequently, $\cA$ is isomorphic to an 
S-ring over $\mZ_{2^{}}\times\mZ_{2^n}$.\medskip

{\bf Case 2: $H$ is dihedral}. In this case all basic sets of $\cA$
other than $\{e\}$ are mixed. Moreover, the S-ring $\cA_H$ is primitive
by the minimality of $H$. Therefore, $\rk(\cA_H)=2$ by 
Theorem~\ref{030414a}. In particular, $\aut(\cA_H)=\sym(H)$. Below,
we will prove that
\qtnl{011014c}
H\le\rad(X)\quad\text{for all}\ X\in\cS(\cA)_{D\setminus H}.
\eqtn
Then the Cayley scheme associated with $\cA$ is isomorphic to
the wreath product of the scheme associated with $\cA_H$ and a
circulant scheme on the right $H$-cosets. Therefore, the group $\aut(\cA)$
contains a subgroup isomorphic to $\sym(H)\wr\mZ_m$, where $m=[D:H]$.
Since the latter group contains a regular subgroup isomorphic to 
$\mZ_{2^{}}\times\mZ_{2^n}$, we are done.\medskip

To complete the proof, we will check statement~\eqref{011014c} in two 
steps: first for rational S-rings, and then in general.

\lmml{280813a} 
Statement~\eqref{011014c} holds whenever the S-ring $\cA$ is rational.
\elmm
\proof  The rationality of $\cA$ implies that it is symmetric, and 
hence commutative. Toward to a contradiction,
suppose that $HX\ne X$ for some basic set $X$ contained in $D\setminus H$. Then 
the product $HX$ is a union of $m>1$ basic sets $X,Y,\ldots$. 
Without loss of generality, we may assume that 
$|X|\leq |Y|\leq\cdots$.\medskip

Since $H$ is an $\cA$-group, it follows from Lemma~\ref{090608a} that
the number $\lambda=|X\cap xH|$ does not depend on the choice of $x\in X$. Therefore,
each $x$ appears $\lambda$ times in the product $\und{H_{}}\,\und{X_{}}$,
i.e.
\qtnl{300514a}
\und{H_{}}\,\und{X_{}} = \lambda(\und{X}+\und{Y}+\cdots).
\eqtn
By the minimality of $X$, this implies that
$|H|\,|X|\geq\lambda m |X|$, and hence $|H|\geq \lambda m$. On the 
other hand, $(H\cap C) X_0 = X_0$ by the rationality of $X$.  Therefore,
the element $\und{X}$  appears  in the product $\und{H}\,\und{X}$ 
at least $|H\cap C|=|H|/2$ times. Thus, $\lambda\geq |H|/2$,
and
$$
|H|\geq \lambda m\ge m|H|/2.
$$
Due to $m>1$, we have $m=2$  and $\lambda=|H|/2$. Consequently, $H_0=H_1=H\cap C$, because
the group $H$ is dihedral. Therefore,
$$
\und{H_{}}\,\und{X_{}} = \und{H_0}(e+s)(\und{X_0}+s\und{X_1}) = 
$$
$$
\und{H_0}\,\und{X_0} +\und{H_0}\,\und{sX_0} +\und{H_0}\,\und{X_1} +\und{H_0}\,\und{sX_1} =
|H_0|\und{X_0} +\und{H_0}\,\und{X_1}+\cdots
$$
By~\eqref{300514a}, all coefficients in the last expression are equal 
to $\lambda=|H_0|$. Therefore, the set $H_0X_1\cap X_0$ must be empty.
It follows that $\und{H_0}\,\und{X_1} = |H_0|\und{H_0X_1}$. However,
this means that $H_0\le\rad(X_1)$. Since also $H_0\le\rad(X_0)$,
we conclude that $H_0\leq\rad(X)$. But $H_0\ne e$, because $H$ is
dihedral. Thus, $\rad(X)$ is non-trivial, and so contains the minimal
$\cA$-group~$H$. But then, $HX = X$. Contradiction.\bull

To complete the proof of~\eqref{011014c}, take a basic set $X$ outside $H$. Then $Y:=\tr(X)$ is also
outside $H$. So by Lemma~\ref{280813a}, we have
$$
\und{Y_0}+\und{sY_1}=\frac{1}{|H|}\und{H}\,\und{Y}=
\frac{1}{|H|}(e+s)\,(\und{H_0Y_0}+\und{H_0Y_1}).
$$ 
Therefore, $|Y_0|=|Y_1|$. On the other hand, for every integer $m$ coprime to $|D|$, we have $|(X^{(m)})_0|=|X_0|$. By
Lemma~\ref{290514b}, this implies that $|(X^{(m)})_1|=|X_1|$. Thus,
$$
k|X_0|=|Y_0|=|Y_1|=k|X_1|,
$$
where $k$ is the number of all distinct sets $X^{(m)}$'s. It follows that $|X_0|=|X_1|$ for each 
basic set $X$ outside $H$.\medskip

Denote by $\rho$ the restriction to $\cA$ of the one-dimensional representation of $D$ that takes $s$ and $c$ to $-1$ and 
$1$, respectively. Then $\rho$ is an irreducible representation of $\cA$ such that $\rho(e)=1$ and $\rho(\und{H^\#})=-1$. Moreover,
for any basic set $X$ outside $H$, we obtain by above that $\rho(\und{X}) =-|X_1|+|X_0| = 0$. In particular, 
$\rho(\und{X^{-1}})=0$. Therefore,
$$
0=\rho(\und{X^{}}\,\und{X}^{-1})=\sum_{Y\in\cS(\cA)} c_{X^{}X^{-1}}^Y\rho(\und{Y})=
$$
$$
c_{X^{}X^{-1}}^e\rho(e)+c_{X^{}X^{-1}}^{H^\#}\rho(H^\#)=
c_{X^{}X^{-1}}^e-c_{X^{}X^{-1}}^{H^\#}.
$$
It follows that $|X|=c_{X^{}X^{-1}}^e=c_{X^{}X^{-1}}^{H^\#}$. Therefore, $HX = X$ for all basis sets
outside $H$, and we are done.\bull

\section{Proof of Theorem~\ref{170814c}}\label{231014a7}

Let $D$ be a dihedral $2$-group and $C$ its cyclic subgroup of index $2$.
Let $\cA$ be an S-ring over the group $D$. Suppose that $r:=\rk(\cA)$ is
at most $5$. For $r=2$, part (1) of Theorem~\ref{170814c} holds trivially. Let $r\ge 3$. 
Then from Theorem~\ref{030414a}, it follows that the S-ring $\cA$ is imprimitive; denote  
by $H$ a minimal non-trivial $\cA$-group. Then $\rk(\cA_H)=2$. Now,
if $r=3$, then $\cA$ is a proper wreath product by Corollary~\ref{021014a}. 
Thus, we can assume that $r=4$ or $r=5$.

\lmml{300514e}
If there is a minimal $\cA$-group $L\ne H$, then statement~(2) holds.
\elmm
\proof By the minimality of the groups $H$ and $L$, we have $H\cap L=e$. Therefore, at least one of them intersects~$C$ trivially. Moreover,
if $H\cap C=L\cap C=e$, then $\grp{HL}$ is an $\cA$-group contained in~$C$, and we replace $H$ by a minimal $\cA$-subgroup in 
$\grp{HL}$. Thus, without loss of generality, we can assume that 
$$
H\cap C\ne e\qaq L\cap C=e.
$$ 
Then $L=\grp{s}$ for some involution 
$s\in D\setminus C$. Moreover, $sHs=H$ by the minimality of $H$. 
Thus, $HL$ is an $\cA$-group and the set $\cS(\cA_{HL})$ contains $4$ elements: $\{e\}$, $H^\#$, $\{s\}$, 
and $sH^\#$. Therefore,
$$
\cA_{HL}=\cA_H\cdot\cA_L.
$$
This implies the required statement for $r=4$, and by Corollary~\ref{021014a} 
also for $r=5$.\bull

By Lemma~\ref{300514e}, we can assume that $H$ is a unique minimal
$\cA$-group. If it is not contained in $C$, then statement (1) holds
by Theorem~\ref{011014a}. Thus, from now, on we also assume that $H\le C$.
Denote by $F$ the union of all basic sets of $\cA$ that are not $C$-mixed.
Clearly, $H\subseteq F$.

\lmml{240913c}
$F$ is an $\cA$-group. Moreover, if $r=r(\cA_F)+1$, then 
$\cA$ is a proper wreath product.
\elmm
\proof The second part of our statement follows from the first one and
Corollary~\ref{021014a}. To prove the first statement, denote by $U$ 
and $V$ the unions of all basic sets of $\cA$ contained 
in $C$ and $D\setminus C$, respectively. We have to prove
that $U\cup V$ is a group. Since $U$ is, obviously, an $\cA$-group,  
without loss of generality, we may assume that $V$ is not empty. Then  
$V=U's$, where $s\in D\setminus C$ is such that $U\cap Us$ is a
subgroup of $D$. It follows that
$$
UU'\subseteq U'\qaq U'U'\subseteq U. 
$$
Since also $U'\subseteq UU'$,  the first inclusion implies that $U'=UU'$. 
Therefore, $U'$ is a union of some $U$-cosets contained in $C$. Since 
the group $C$ is a cyclic $2$-group, this together with the 
second inclusion implies $U'=U$. Thus, the set 
$U\cup V=U\cup Us$ is a group.\bull

Suppose first that $F_0=C$. Then from the definition of $F$, it follows that $C$ is an 
$\cA$-group. By Corollary~\ref{021014a}, we can assume that $r_C=\rk(\cA_C)$
is not equal to $r-1$. Since, obviously, $r_C\ge 2$,
we have
$$
(r,r_C)=(4,2),\ (5,2)\ \text{or}\ (5,3).
$$
In the former two cases, we are done by statement~(1) of Corollary~\ref{120614a},
whereas in the third one by statement~(2). Thus, in
what follows, we can assume that 
$$
F_0<C\qaq H=F\ \text{or}\  r_F=3,
$$
where $r_F=\rk(\cA_F)$. In particular, there are two 
or three basic sets outside $F$ (notice, that they are $C$-mixed). 

\lmml{041014a}
Let $X\in\cS(\cA)_{D\setminus F}$. Suppose that $X$ is rational or 
$[V:H]\ge 4$, where $V=\grp{X_0}$. Then $H\le\rad(X)$. 
\elmm
\proof It suffices to verify that $H\le\rad(X_0)$. Indeed, then
the coefficient at $\und{X}$ in $\und{H}\,\und{X}$ is at least $|H|$.
Since it can not be larger than $|H|$, we are done.\medskip

Suppose first that $X$ is rational. Then by statement~(2) of Lemma~\ref{240913a} the
set $X_0$ is rational. Since $X_0\subseteq C\setminus H$, we conclude that
$H\le\rad(X_0)$, as required.\medskip

Let now $[V:H]\ge 4$. Then $V\cong\mZ_{2^k}$ for some $k\ge 2$. Since $r\le 5$, there are at most two basic 
rationally conjugate to $X$. Therefore, the stabilizer of $X_0$ in the group $(\mZ_{2^k})^*$, 
has index at most $2$ in it. It follows that this stabilizer 
contains the subgroup of all elements $x\mapsto x^{1+4m}$, $x\in \mZ_{2^k}$,
with $m\in \mZ_{2^k}$. By statement~(2) of Lemma~\ref{211113a}, this implies that 
$\rad(X_0)\ge V^4\ge H$.\bull

From Lemma~\ref{041014a}, it follows that if all basic sets outside $F$ are
rational, then $\cA$ is a proper wreath product. Indeed, this is obvious when $H=F$.
If $H\ne F$, this is also true, because then $F\setminus H$ is a basic set, the radical 
of which equals $H$. Thus, we can assume that two of basic sets outside $F$, say $X$ and $Y$, 
are rationally conjugate, and the third one (if exists) is rational. The rest of the proof is divided 
into four cases below.\medskip

{\bf Case 1: $F=H$ and $r=4$.}
Using the computer package COCO, \cite{PRZ} we found exactly five and three
S-rings of rank~$4$ over the  groups $D_8$ and $D_{16}$, respectively. 
In both cases, only two of them have a unique minimal $\cA$-group 
contained in $C$ and they are proper wreath products. Thus, in what follows, we assume 
that $|D|\ge 32$.\medskip

In our case,  the non-trivial basic sets of $\cA$ are $X$, $Y$ and $Z=H^\#$. It is easily seen that
the hypothesis of Lemma~\ref{290514b} is satisfied. Since $X$ and $Y$
are rationally conjugate, there is an algebraic isomorphism of $\cA$ that takes $X$
to $Y$, and $Y$ to $X$. Therefore,
\qtnl{081014a}
\und{H}\,\und{X}=a(\und{X}+\und{Y}),
\eqtn
where $a=|H|/2$. Consequently, $c_{ZX}^X=a-1$. On the other hand, since $X$ and $Z$ are 
symmetric, we have $C_{XX}^Z=\frac{|X|}{|Z|}C_{ZX}^X$.  It follows that $|Z|=2a-1$ divides 
$$
|X|c_{ZX}^X=\frac{(d-2a)(a-1)}{2},
$$
where $d=|D|$. However, by Lemma~\ref{041014a}, without loss of generality,
we can assume that $|C:H|=2$. Therefore, $2a=|H|=d/4$. It follows that
$2a-1$ divides $3a(a-1)$. But $a$ being a power of $2$, must be coprime to $2a-1$.
Consequently, $2a-1$ divides $3a-3$. Since this is possible only for $a\le 2$, i.e.
when $d\le 16$, we are done.\medskip

{\bf Case 2: $F=H$ and $r=5$.} Denote by $Z$ the basic set in $\cS(\cA)_{D\setminus H}$
other than $X$ and $Y$. It is easily seen that the hypothesis of Lemma~\ref{290514b} 
is satisfied. Since $X$ and $Y$ are rationally conjugate, there is an algebraic isomorphism of $\cA$ 
that takes $X$ to $Y$, $Y$ to $X$ and leaves $Z$ fixed. Therefore, the rational closure of
$\cA$ is of rank~$4$. So by Lemma~\ref{041014a}, it is the
wreath product $\cA_H\wr\cB$, where $\cB$ is isomorphic to the rational closure of $\cA_{D/H}$.
Therefore, $\rk(\cB)=3$, and hence there exists  a non-trivial $\cB$-group. Denote by $U$
its preimage in $\cA$. Then, obviously, $H<U<D$.\medskip

Since $X$ and $Y$ are rationally conjugate, we have $X\cup Y\subseteq U$ or 
$X\cup Y\subseteq D\setminus U$. However, $\rk(\cA)=r=5$. Therefore, $Z=D\setminus U$ in the former 
case, and $Z=U\setminus H$ in the latter one. In any case, $H\le\rad Z$. By Lemma~\ref{041014a},
this implies that if $Z=U\setminus H$, then 
\qtnl{140215a}
\rad(X)=\rad(Y)\ge H,
\eqtn
and $\cA$ is a proper wreath product. Let now, $Z=D\setminus U$. Then $\rk(\cA_U)=4$, and $\cA_U$ is
the wreath product $\cA_H\wr\cA_{U/H}$. Therefore, again \eqref{140215a} holds, and $\cA$ is a proper wreath product.\medskip

{\bf Case 3: $C\not\ge F>H$.} In this case, $F=H\cup Hs$. So by Lemma~\ref{041014a},
the sets $X_0$ and $Y_0$ are orbits  of an index $2$ subgroup of $\aut(C)$, unless
$\cA$ is a proper wreath product.
This implies that the group $\rad(X_0)=\rad(Y_0)$ has index $2$ in $H$. Therefore,
equality~\eqref{081014a} holds with $a=|H|/2$. Exactly as in Case~2, we conclude
that $2a-1$ divides 
$$
|X|c_{ZX}^X=\frac{(d-4a)(a-1)}{2},
$$
where $Z=H^\#$ and $d=|D|=4|H|=8a$. Thus, $2a-1$ divides $2a(a-1)$. Contradiction.\medskip

{\bf Case 4: $C\ge F>H$.} In this case, $F=F_0<C$ by the above assumption. Therefore,
$X_0$ (and also $Y_0$) contains a generator of $C$. It follows that
$$
[\grp{X_0}:H]=[C:H]\ge [C:F][F:H]\ge 4.
$$
By Lemma~\ref{041014a}, this implies that \eqref{140215a} holds. Since also
$F\setminus H$ is the basic set and $H=\rad(F\setminus H)$, the group $H$ is contained
in the radical of every basic set outside~$H$. Thus, $\cA$ is a proper wreath product.\bull


\begin{thebibliography}{99}

\bibitem{BC05}
R.~A.~Bailey, P.~J.~Cameron, {\em Crested products of association schemes},
J. London Math. Soc., {\bf 72} (2005), 2, 1-–24.

\bibitem{BJL}
T.~Beth, D.~Jungnickel, H.~Lenz, {\em Design Theory}, 2nd edition, Cambridge University 
Press, Cambridge, 1999.

\bibitem{C72}
P.~J.~Cameron, {\em On groups with several doubly-transitive permutation representations},
Math. Z., {\bf 128} (1972), 1--14.

\bibitem{DS94}
J.~A.~Davis, K.~Smith, {\em A construction of difference sets in high exponent 2-groups using
representation theory}, Journal Algebraic Combin., {\bf 3} (1994), 137--151.

\bibitem{EP02}
S.~Evdokimov, I.~Ponomarenko, {\em Characterization of cyclotomic schemes and normal Schur rings over a cyclic group}, St.
Petersburg Math. J.,  {\bf 14}  (2003),  no. 2, 189--221.

\bibitem{EP05}
S.~Evdokimov, I.~Ponomarenko, {\em A new look at the Burnside-Schur
theorem}, Bulletin of the London Mathematical Society, {\bf 37} (2005), 535–-546.

\bibitem{EP09}
S.~Evdokimov, I.~Ponomarenko, {\em Permutation group approach to association schemes},
European Journal Combin., {\bf 30} (2009), 6, 1456--1476.

\bibitem{EvdP09}
S.~Evdokimov, I.~Ponomarenko, {\em Schur rings over a product of Galois rings},
Beitr. Algebra Geom., {\bf 55} (2014), 1, 105--138. 

\bibitem{EP12}
S.~Evdokimov, I.~Ponomarenko, {\em Schurity of S-rings over a cyclic group and generalized 
wreath product of permutation groups}, Algebra and Analysis, {\bf 24} (2012), 3, 84--127.

\bibitem{EKP} 
S.~Evdokimov, I.~Kov\'acs, I.~Ponomarenko, {\em Characterization 
of cyclic Schur groups}, Algebra and Analysis, {\bf 25} (2013), 5, 61--85. 

\bibitem{EKP14}
S.~Evdokimov, I.~Kov\'acs, I.~Ponomarenko, {\em On schurity of finite abelian groups},
{\tt arXiv:1309.0989 [math.GR]} (2013), 1--20 (accepted to Communication in Algebra).

\bibitem{HM}
M.~Hirasaka, M.~Muzychuk, {\em Association schemes generated by a non-symmetric relation 
of valency 2}, Discrete Math., {\bf 244} (2002), 109--135.

\bibitem{Jo02}
G.~A.~Jones, {\em Cyclic regular subgroups of primitive permutation groups},
Journal Group Theory, {\bf 5} (2002), 403–-407.

\bibitem{KK}
L.~A.~Kaluzhnin, M.~Kh.~Klin, {\em Certain numerical invariants of permutation groups},
Latv. Mat. Ezheg., {\bf 18} (1976), 81--99 (in Russian).

\bibitem{K85}
M.~Kh.~Klin, {\em The axiomatics of cellular rings}, Investigations in the
algebraic theory of combinatorial objects, 6--32, Vsesoyuz.
Nauchno-Issled. Inst. Sistem. Issled., Moscow, 1985 (in Russian).

\bibitem{K37}
R.~Kochend{\"o}rffer, \emph{Untersuchungen {\"u}ber eine Vermutung 
von W. Burnside}, Schr. Math. Semin. u. Inst. Angew. Math. Univ. 
Berlin, {\bf 3} (1937), 155--180.


\bibitem{KMM11} 
I.~Kov{\' a}cs, D.~Maru{\u s}i{\u c}, M.~Muzychuk, {\em On dihedrants admitting 
arc-regular group actions}, Journal Algebraic Combin., {\bf 35} (2011), 409--426. 

\bibitem{LeungMa90} 
K.~H.~Leung, S.~L.~Ma, {\em The structure of Schur rings over cyclic groups},
Journal Pure Appl. Algebra, {\bf 66} (1990), 287--302.

\bibitem{LM98}
K.~H.~Leung, S.~H.~Man, {\em On Schur Rings over Cyclic Groups},
Israel J. Math., {\bf 106} (1998), 251--267.

\bibitem{M87}
M.~Muzychuk, {\em $V$-rings of permutation groups that have invariant metric}, 
Thesis, (1987) (in Russian).

\bibitem{MP09}
 M.~Muzychuk, I.~Ponomarenko, {\em Schur rings}, European Journal Combin.,
{\bf 30} (2009), 6, 1526--1539.

\bibitem{MP}
M.~Muzychuk, I.~Ponomarenko, {\em On quasi-thin association schemes},
Journal of Algebra, {\bf 351} (2012), 467--489.

\bibitem{PRZ}
Ch.~Pech, S.~Reichard, M.~Ziv-Av, {\em The COCO share package for
{\sf GAP}}, https://github.com/MatanZ/coco-ii.
 
\bibitem{P13}
I. ~N.~Ponomarenko, {\em Bases of Schurian antisymmetric coherent configurations and isomorphism 
test for Schurian tournaments}, J. Math. Sci. (N. Y.), {\bf 192} (2013), no. 3, 316-–338.

\bibitem{PV}
I.~Ponomarenko,  A.~Vasil'ev, {\em On non-abelian Schur groups},
Journal of Algebra and Its Applications, {\bf 13} (2014), 8, 1450055-1--1450055-22.

\bibitem{Poe74}
R.~P\"oschel, {\em Untersuchungen von S-ringen insbesondere im gruppenring von
$p$-gruppen}, Math. Nachr.,  {\bf 60} (1974), 1--27.

\bibitem{P95}
A.~Pott, {\em Finite geometry and character theory}, Berlin: Springer-Verlag, 1995.

\bibitem{S33}
I.~Schur, {\em Zur Theorie der einfach transitiven Permutationgruppen},
S.-B. Preus Akad. Wiss. Phys.-Math. Kl., (1933), 598--623.


\bibitem{W49}
H.~Wielandt, \emph{Zur Theorie der einfach transitiven Permutationsgruppen. II},
Mathematische Zeitschrift, {\bf 52} (1949), 384--393.

\bibitem{Wie64}
H.~Wielandt, {\em Finite permutation groups}, Academic Press,
New York - London, 1964.

\end{thebibliography}
\end{document}